\documentclass[12pt,letterpaper]{article}
\usepackage[T1]{fontenc}
\usepackage{mathtools,amsthm,amssymb}
\usepackage{txfonts}
\usepackage{microtype}
\usepackage[letterpaper,textwidth=160mm,textheight=235mm,centering]{geometry}
\usepackage{indentfirst,needspace,cite}
\usepackage[colorlinks=true,linkcolor=blue,urlcolor=blue,citecolor=blue]{hyperref}

\newtheorem{thm}{Theorem}[section]
\newtheorem{lem}{Lemma}[section]
\newtheorem{prob}{Problem}[section]
\newtheorem{claim}{Claim}

\newcommand{\ceil}[1]{\left\lceil #1\right\rceil}
\newcommand{\floor}[1]{\left\lfloor #1\right\rfloor}
\newcounter{casegroup}
\newcounter{case}[casegroup]
\newcounter{subcase}[case]
\renewcommand{\thesubcase}{\thecase.\arabic{subcase}}

\newcommand{\resetcases}{\stepcounter{casegroup}}
\newenvironment{case}[2][1.5mm]{%
  \vspace{#1}%
  \noindent\refstepcounter{case}\textbf{Case \thecase.} #2%
  \space\vspace{1mm}%
}{}
\newenvironment{subcase}[2][1.5mm]{%
  \vspace{#1}%
  \noindent\refstepcounter{subcase}\textbf{Subcase \thesubcase.} #2%
  \space\vspace{1mm}%
}{}

\makeatletter
\renewcommand{\@maketitle}{%
  \newpage\null\vskip 1em
  \begin{center}
    {\Large\@title\par}%
    \ifx\@author\@empty\else
      \vskip 1em{\normalsize\begin{tabular}[t]{c}\@author\end{tabular}\par}%
    \fi
    \ifx\@date\@empty\else\vskip .75em{\normalsize\@date\par}\fi
  \end{center}\par\vskip 1em}
\makeatother

\begin{document}

\title{Extremal problems for graphs with bounded $k$-th binding number\footnote{Supported by the National Natural Science Foundation of China
{(No. 12371361)}, and Distinguished Youth Foundation of Henan Province {(No. 242300421045)}.}}
\author{{\bfseries Ao Fan, Ruifang Liu}\thanks{Corresponding author.
E-mail addresses: rfliu@zzu.edu.cn (R. Liu), fanaozzu@163.com (A. Fan).}\\
{\footnotesize School of Mathematics and Statistics, Zhengzhou University, Zhengzhou, Henan 450001, China}}
\date{}
\maketitle
{\flushleft\large\bf Abstract}
The binding number of a graph, introduced by Woodall [J. Combin. Theory, Ser. B, 1973], is a central topic in structural and extremal graph theory.
Recently, Chen et al. introduced the $k$-th binding number $\beta^k(G)$. Anderson [J. Combin. Theory, Ser. B, 1971] proved that every graph $G$ of even order with $\beta^1(G)\ge\frac43$ contains a perfect matching, and that the bound $\frac43$ is best possible. For $k\ge2$, Chen et al. proved that $\beta^k(G)\ge1$ guarantees a $k$-factor whenever $kn$ is even, and that the bound $1$ is best possible. The above facts reveal a fundamental distinction between the higher-order binding number and the classical binding number. This distinction makes extremal problems for graphs with bounded $k$-th binding number particularly intriguing.

For positive integers $n,$ $k$ and $r$, we in this paper determine the maximum size and spectral radius among all graphs $G$ of order $n$ satisfying  $\beta^k(G)<r$, and characterize the corresponding extremal graphs, respectively. In particular, when $k=1$, our results imply the corresponding
extremal results of Fan and Lin [Electron. J. Combin., 2024] and Liu, Chen and Fan [Electron. J. Combin., 2026]. Our proofs develop a structural framework that relates higher-order neighborhood inequalities to the structure of the complements of graphs. This framework may provide an approach to Tur\'{a}n-type problem involving the $k$-th binding number, and the approach may be of independent interest.

\begin{flushleft}
\end{flushleft}

\section{Introduction}

Let $G$ be a simple graph with vertex set $V(G)$ and edge set $E(G)$.
The order and size of $G$ are denoted by $|V(G)|=n$ and $|E(G)|=e(G)$,
respectively. We write $K_n$, $P_n$ and $C_n$ for the complete
graph, path and cycle of order $n$, respectively. For a vertex $v\in V(G)$, let $\Lambda_G(v)$ denote its
neighborhood.
For a vertex subset $S\subseteq V(G)$, let $G[S]$ be the subgraph
of $G$ induced by $S$, and write $\Lambda_G(S)=\bigcup_{v\in S}\Lambda_G(v)$.
The complement of $G$ is denoted by $\overline G$.
Let $G_1$ and $G_2$ be two vertex-disjoint graphs. Their disjoint union
is denoted by $G_1\cup G_2$, and their join $G_1\vee G_2$ is obtained
from $G_1\cup G_2$ by adding all edges between $V(G_1)$ and $V(G_2)$.
The largest eigenvalue of the adjacency matrix $A(G)$ is called
the \emph{spectral radius} of $G$, and is denoted by $\rho(G)$.
For undefined terminology and notation, one can refer to
\cite{Bondy2008,Godsil2001}.

Woodall \cite{Woodall1973} introduced the \emph{binding number}
$b(G)$ of a graph $G$, defined as
$$
b(G)=\min\left\{\frac{|\Lambda_G(S)|}{|S|}:
\emptyset\ne S\subseteq V(G),\ \Lambda_G(S)\ne V(G)\right\}.
$$
A set attaining this minimum is called a \emph{binding set}.
Recently, Chen et al. \cite{Chen2026} introduced the
\emph{$k$-th binding number}.
For a positive integer $k$ and a vertex subset $S\subseteq V(G)$, let
$
d_S(v)=|\Lambda_G(v)\cap S|,$
$\Lambda_G^k(S)=\{v\in V(G):d_S(v)\ge k\}.
$
The set $\Lambda_G^k(S)$ is called the \emph{$k$-th neighborhood} of $S$.
For $n\ge k$, the $k$-th binding number of $G$ is defined as
\begin{equation*}\label{eq:def}
\beta^k(G)=\min\left\{\frac{|\Lambda_G^k(S)|}{|S|}:
S\subseteq V(G),\ |S|\ge k,\ \Lambda_G^k(S)\ne V(G)\right\}.
\end{equation*}
Following \cite{Chen2026}, we put $\beta^k(G)=0$ when $n<k$.
For convenience, we write $\Lambda^k(S)$ in place of $\Lambda_G^k(S)$
whenever the graph $G$ is clear from the context.
For $k=1$, we have $\Lambda_G^1(S)=\Lambda_G(S)$, and hence
$\beta^1(G)=b(G).$

As a fundamental graph invariant, the binding number possesses wide-ranging applications in numerous fields, such as graph theory, network science, quantum sensing, and information processing. A central problem is to understand its relationship with basic structural parameters, including
the degree sequence \cite{Bauer2011}, minimum degree \cite{Bauer2012} and connectivity \cite{Woodall1973}.
Lower bounds on the binding number $b(G)$
also give sufficient conditions for several structural properties.
Anderson \cite{Anderson1971} gave a condition ensuring a perfect
matching in a graph of even order.
Woodall \cite{Woodall1973,Woodall1978} established sufficient
conditions for Hamiltonicity, and Shi \cite{Shi1987} obtained
a corresponding result for pancyclicity.
Katerinis and Woodall \cite{Katerinis1987} studied $k$-factors,
and Kano and Tokushige \cite{Kano1992} considered $f$-factors.
Chen \cite{Chen1995} and Robertshaw and Woodall
\cite{Robertshaw2002} obtained sufficient conditions for
$k$-extendibility, ensuring that every matching of size $k$
extends to a perfect matching.
The binding number is also related to clique containment.
Shi \cite{Shi1985} gave a sufficient condition for the existence
of a triangle, and Lyle and Goddard \cite{Lyle2009} improved
upper bounds on the binding number of graphs with no clique
of a given order.

Anderson \cite{Anderson1971} proved that every graph $G$
of even order with $\beta^1(G)=b(G)\ge\frac43$ contains a $1$-factor,
and that the lower bound $\frac43$ is best possible.
For $k=2$, Bar\'at et al. \cite{Barat2024} proved that
$\beta^2(G)\ge1$ guarantees a $2$-factor in $G$.
The lower bound $1$ is best possible.
For general $k\ge2$, Chen et al. \cite{Chen2026} proved that,
when $k|V(G)|$ is even, the condition $\beta^k(G)\ge1$ guarantees
a $k$-factor in $G$. The lower bound $1$ is best possible.
Moreover, Chen et al. \cite{Chen2026} proved that any graph $G$ with
$\beta^k(G)\ge1$ contains at least $k-1$ disjoint perfect or near-perfect matchings.
The above facts reveal a fundamental distinction between the higher-order binding number and the classical binding number.

Fan and Lin \cite{Fan2024} determined the maximum size or
spectral radius among connected graphs of order $n$
satisfying $b(G)<r$, where $r\ge1$ is an integer,
and characterized the extremal graphs.
Liu, Chen and Fan \cite{Liu2026} determined the maximum size
or spectral radius for graphs without isolated vertices
satisfying $b(G)<\frac{1}{r}$.
Furthermore, they also determined the extremal graphs maximizing the
size or the spectral radius among all bipartite graphs with
$b(G)<\frac{1}{r}$.
The fundamental differences between the classical binding number and its higher-order counterparts motivate us to investigate analogues of these extremal results for the $k$-th binding number. We therefore pose the following problem.

\begin{prob}\label{prob:main}
For positive integers $k$ and $r$, what is the maximum size or
spectral radius among graphs $G$ of order $n$ without isolated vertices
satisfying $\beta^k(G)<r$?
\end{prob}

In this paper, we completely solve Problem~\ref{prob:main} and
characterize all extremal graphs.
Throughout the paper, let $n,k,r$ be three positive integers.
Let $\mathcal G_{n,k,r}$ be the family of graphs $G$ of order $n$
such that $\delta(G)\ge1$ and $\beta^k(G)<r$.
Note that, for $n\ge k$, every vertex subset $S$ of size $k$
satisfies $S\cap\Lambda^k(S)=\emptyset$. Hence
$\beta^k(G)\le\frac{n-k}{k}$ and $\beta^k(K_n)=\frac{n-k}{k}.$
Together with $\beta^k(G)=0$ for $n<k$, this implies that every
graph of order $2\le n<(r+1)k$ satisfies $\beta^k(G)<r$.
It is easy to see that $K_n$ is the unique graph attaining the
maximum size or spectral radius in $\mathcal G_{n,k,r}$ for $2\le n<(r+1)k$. Therefore, it suffices to consider $n\ge(r+1)k$.

Since $e(G)=\binom n2-e(\overline G)$, maximizing the size of
$G$ is equivalent to minimizing the number of missing edges.
Here and below, $K_n-F$ denotes the spanning subgraph of $K_n$
obtained by deleting the edges of a copy of $F$. For $1\le t\le n-2$, let $H_{n,t}=K_{n-t-1}\vee(K_t\cup K_1).$ For $(k,r)\ne(1,1)$ and $n\ge(r+1)k$, define
\begin{equation*}\label{eq:p}
p=p(n,k,r)=
\begin{cases}
1,&(r+1)k\leq n\le r(k+1)+1,\\[1mm]
\min\{n-(r+1)k+1,\ceil{\frac{n}{r}}-k\},&n\ge r(k+1)+2.
\end{cases}
\end{equation*}
Let $\mathcal S_{k,t}$ be the family of star forests with $t$ edges,
no isolated vertices and at most $k$ components, and let $\mathcal F_{n,k,r}$ be the
family of graphs $F$ with $p$ edges and no isolated vertices
satisfying at least one of the following conditions.\\
\textup{(a)} $F\cong K_{1,p}$.\\
\textup{(b)} $F\in\mathcal S_{k,p},$ if $n=(r+1)k+p-1$.\\
\textup{(c)} $F\cong2K_2$, if $k\ge3,$ $n\le r(k+1)+3$ and $p=2$.\\
\textup{(d)} $F\cong P_3\cup K_2$, if $k\ge4,$ $n\le r(k+1)+4$ and $p=3$.\\
\textup{(e)} $F\cong3K_2$, if $k\ge5$, $n\le r(k+1)+5$ and $p=3$.\\
\textup{(f)} $F\cong P_4$, if $k\ge2$, $n\le r(k+2)+1$ and $p=3$.\\
\textup{(g)} $F\cong K_3$, if $n\le r(k+2)+2$ and $p=3$.\\
\textup{(h)} $F\cong C_4$, if $k\ge2$, $n\le r(k+2)+3$ and $p=4$.\\

Recall that $\mathcal G_{n,k,r}$ is the family of graphs $G$ of order $n$ such that $\delta(G)\ge1$ and $\beta^k(G)<r$. 

\begin{thm}\label{edge:main}
Let $k,r\ge1$ and $n\ge(r+1)k$ be integers with $(k,r)\ne(1,1)$,
and let $G\in\mathcal G_{n,k,r}$.\\
(i) If $(n,k,r)=(11,3,2)$, then $e(G)\le53$ with equality if and
only if $G\cong K_{11}-2K_2$.\\
(ii) If $(n,k,r)\ne(11,3,2)$, then
$e(G)\le\binom n2-p$
with equality if and only if $G\cong K_n-F$ for some
$F\in\mathcal F_{n,k,r}$.
\end{thm}

We next consider the remaining case $k=r=1$.
The family $\mathcal G_{2,1,1}$ is empty, since $K_2$ is the
only graph of order two without isolated vertices and
$b(K_2)=1$. We only consider $n\ge3$.

\Needspace{10\baselineskip}
\begin{thm}\label{edge:first}
Let $G$ be a graph of order $n\ge3$ without isolated vertices
and with $\beta^1(G)=b(G)<1$, and set
$q=\floor{\frac{n-1}{2}}$.\\
(i) If $n\in\{5,6,7,8,9,11\}$, then
$e(G)\le\binom q2+q(n-q)$ with equality if and only if
$G\cong K_q\vee(n-q)K_1$.\\
(ii) If $n=10$, then $e(G)\le30$ with equality if and
only if $G\cong K_1\vee(K_7\cup2K_1)$ or
$G\cong K_4\vee6K_1$.\\
(iii) If $n=13$, then $e(G)\le57$ with equality if and
only if $G\cong K_1\vee(K_{10}\cup2K_1)$ or
$G\cong K_6\vee7K_1$.\\
(iv) If $n\in\{3,4,12\}$ or $n\ge14$, then
$e(G)\le\binom{n-2}{2}+2$ with equality if and only if
$G\cong K_1\vee(K_{n-3}\cup2K_1)$.
\end{thm}

\Needspace{9\baselineskip}

For $k=1$, Theorems~\ref{edge:main} and \ref{edge:first}
extend Fan and Lin's size result to graphs of general order $n$ without isolated
vertices.

\begin{thm}[Fan and Lin {\cite{Fan2024}}]\label{thm:fanlin-edge}
Let $r$ be a positive integer, and let $G$ be a connected graph
of order $n$ with $b(G)<r$.\\
(i) If $r=1$ and $n\ge14$, then
$e(G)\le\binom{n-2}{2}+2$ with equality if and only if
$G\cong K_1\vee(K_{n-3}\cup2K_1)$.\\
(ii) If $r\ge2$ and $n\ge3r+5$, then
$e(G)\le\binom n2-\ceil{\frac{n-r}{r}}$
with equality if and only if
$G\cong K_{n-\ceil{\frac{n-r}{r}}-1}\vee(K_{\ceil{\frac{n-r}{r}}}\cup K_1)$.
\end{thm}

We now turn to the maximum spectral radius of graphs in $\mathcal G_{n,k,r}$.

\begin{thm}\label{thm:main}
Let $k\ge1$, $r\ge1$ and $n\ge(r+1)k$ be integers with
$(k,r)\ne(1,1)$, and let $G\in\mathcal G_{n,k,r}$. \\
(i) If $(n,k,r)=(11,3,2)$, then
$\rho(G)\le\rho(K_{11}-2K_2)$ with equality if and only if
$G\cong K_{11}-2K_2$.\\
(ii) If $(n,k,r)\ne(11,3,2)$, then
$\rho(G)\le\rho(H_{n,p})$ with equality if and only if
$G\cong H_{n,p}$.
\end{thm}

We next determine the maximum spectral radius in the remaining case $k=r=1$.

\begin{thm}\label{thm:first}
Let $G$ be a graph of order $n\ge3$ without isolated vertices
and with $\beta^1(G)=b(G)<1$. \\
(i) If $n=5$, then $\rho(G)\le\rho(K_2\vee3K_1)$
with equality if and only if $G\cong K_2\vee3K_1$.\\
(ii) If $n=6$, then $\rho(G)\le\rho(K_2\vee4K_1)$
with equality if and only if $G\cong K_2\vee4K_1$.\\
(iii) If $n=7$, then $\rho(G)\le\rho(K_3\vee4K_1)$
with equality if and only if $G\cong K_3\vee4K_1$.\\
(iv) If $n=9$, then $\rho(G)\le\rho(K_4\vee5K_1)$
with equality if and only if $G\cong K_4\vee5K_1$.\\
(v) If $n\in\{3,4,8\}$ or $n\ge10$, then
$\rho(G)\le\rho(K_1\vee(K_{n-3}\cup2K_1))$ with equality
if and only if $G\cong K_1\vee(K_{n-3}\cup2K_1)$.
\end{thm}

Let $k,r\ge1$ be integers, and suppose that
$\mathcal G_{n,k,r}\ne\emptyset$.
If $(k,r)\ne(1,1)$, then the unique graph with the maximum
spectral radius in $\mathcal G_{n,k,r}$ also has the maximum size.
In particular, $K_{11}-2K_2$ attains the maximum size and
spectral radius for $(n,k,r)=(11,3,2)$.
If $k=r=1$, then the spectral extremal graph has the maximum
size if and only if $n\notin\{8,11\}$.
Theorems~\ref{thm:main} and \ref{thm:first}
extend Fan and Lin's spectral result to graphs of general order $n$ without isolated
vertices.
\begin{thm}[Fan and Lin {\cite{Fan2024}}]\label{thm:fanlin-spectral}
Let $r$ be a positive integer, and let $G$ be a connected graph
of order $n$ with $b(G)<r$. \\
(i) If $r=1$ and $n\ge14$, then
$\rho(G)\le\rho(K_1\vee(K_{n-3}\cup2K_1))$ with equality if and
only if $G\cong K_1\vee(K_{n-3}\cup2K_1)$.\\
(ii) If $r\ge2$ and
$n\ge\max\{\frac{r^3+5r^2+r}{2},5r^2+3\}$, then
$\rho(G)\le\rho(K_{n-\ceil{\frac{n-r}{r}}-1}\vee(K_{\ceil{\frac{n-r}{r}}}\cup K_1))$
with equality if and only if
$G\cong K_{n-\ceil{\frac{n-r}{r}}-1}\vee(K_{\ceil{\frac{n-r}{r}}}\cup K_1)$.
\end{thm}

\section{Preliminaries}\label{sec:combinatorial}

In this section, we establish several structural lemmas that will be used in
the proofs of our size and spectral extremal results.
\begin{lem}\label{lem:starbeta}
Let $k\ge1$, $n\ge2k$ and $1\le t\le n-2$. Then
$$
\beta^k(H_{n,t})=\min\left\{
\frac{\max\{0,n-k-t\}}{k},
\frac{n-1-\min\{t,k\}}{k+1},
\frac{n-1}{\min\{n,k+t\}}\right\}.
$$
\end{lem}
\begin{proof}
Let $v$ be a vertex of degree $n-t-1$ in $H_{n,t}$, let $U$
be the set of its $t$ nonneighbors, and put $W=\Lambda_{H_{n,t}}(v)$.
Then $|W|=n-t-1$, and $U\cup W$
induces a complete graph.
Let $S\subseteq V(H_{n,t})$ satisfy $|S|\ge k$ and
$\Lambda^k(S)\ne V(H_{n,t})$.

\resetcases
\begin{case}{$v\in S$.}\label{case:starbeta-in}

If $|S|=k$, then $\Lambda^k(S)=W\setminus S$, and hence
$\frac{|\Lambda^k(S)|}{|S|}=\frac{|W|-|S\cap W|}{k}$.
Thus, the ratio is minimized when $|S\cap W|=\min\{k-1,|W|\}$.
This is attained by choosing $\min\{k-1,|W|\}$ vertices from $W$
and the remaining $k-1-\min\{k-1,|W|\}$ vertices from $U$,
together with $v$.
Therefore, the minimum ratio in this case is
$$
\frac{|W|-\min\{k-1,|W|\}}{k}
=\frac{\max\{0,n-k-t\}}{k}.
$$
If $|S|=k+1$ and $S\cap U\ne\emptyset$, then
$\Lambda^k(S)=W\cup(U\setminus S)$
and $\frac{|\Lambda^k(S)|}{|S|}=\frac{n-1-|S\cap U|}{k+1}$.
Thus, the ratio is minimized by choosing $\min\{t,k\}$ vertices
from $U$ and $k-\min\{t,k\}$ vertices from $W$, together with $v$.
Therefore, the minimum ratio in this case is
$$\frac{n-1-\min\{t,k\}}{k+1}.$$
If $|S|=k+1$ and $S\cap U=\emptyset$, then every vertex of
$H_{n,t}$ has at least $k$ neighbors in $S$.
Hence $\Lambda^k(S)=V(H_{n,t})$, a contradiction.

Now suppose that $|S|\ge k+2$.
Every vertex other than $v$ has at least $|S|-2\ge k$ neighbors
in $S$ and hence belongs to $\Lambda^k(S)$.
Since $\Lambda^k(S)\ne V(H_{n,t})$, we must have
$v\notin\Lambda^k(S)$.
Therefore, $\Lambda^k(S)=V(H_{n,t})\setminus\{v\}$,
$|S\cap W|=d_S(v)\le k-1$.
Since $v\in S$, $|S\cap U|\le t$ and $|S\cap W|\le k-1$,
we have $|S|=1+|S\cap U|+|S\cap W|\le1+t+(k-1)=k+t.$
Therefore, $|S|\le\min\{n,k+t\}$.

If $t=1$, then $|U|=1$. Recall that $v\notin\Lambda^k(S)$,
so $|S\cap W|=d_S(v)\le k-1$.
Since $v\in S$, we obtain
$
|S|=1+|S\cap U|+|S\cap W|\le1+1+(k-1)=k+1,
$
contrary to $|S|\ge k+2$.
Moreover, for $t=1$, we have
$$
\begin{aligned}
\frac{n-1}{\min\{n,k+t\}}
=\frac{n-1}{\min\{n,k+1\}}=\frac{n-1}{k+1}
>\frac{n-2}{k+1}
=\frac{n-1-\min\{t,k\}}{k+1}.
\end{aligned}
$$
Hence the third term is strictly larger than the second and does not
affect the stated minimum.

Next we consider $t\ge2$.
Choose $W'\subseteq W$ with $|W'|=\min\{k-1,|W|\}$,
and put $S=\{v\}\cup U\cup W'$.
Note that $n\ge2k\ge k+2$ for $k\ge2$.
For $k=1$, we have $n\ge t+2\ge4$ since $t\ge2$.
Then
$|S|=1+t+\min\{k-1,n-t-1\}=\min\{n,k+t\}\ge k+2.$
Every vertex other than $v$ has at least $|S|-2\ge k$ neighbors
in $S$, whereas $d_S(v)=|W'|\le k-1$.
Hence $\Lambda^k(S)=V(H_{n,t})\setminus\{v\}\ne V(H_{n,t})$.
Thus, $S$ satisfies the conditions in the definition of $\beta^k(H_{n,t})$ and
attains the lower bound
$$
\frac{|\Lambda^k(S)|}{|S|}=\frac{n-1}{\min\{n,k+t\}}.
$$

\end{case}
\begin{case}[1mm]{$v\notin S$.}\label{case:starbeta-out}

If $|S|=k$, then every vertex of $(U\cup W)\setminus S$
is adjacent to all vertices of $S$, since $U\cup W$
induces a complete graph $K_{n-1}$. Note that $t\ge1$ and $n\ge k+1$. Then
$$
\frac{|\Lambda^k(S)|}{|S|}
\ge\frac{n-k-1}{k}
\ge\frac{\max\{0,n-k-t\}}{k}.
$$

Now suppose that $|S|\ge k+1$.
Every vertex other than $v$ has at least
$|S|-1\ge k$ neighbors in $S$ and hence belongs to $\Lambda^k(S)$.
Since $\Lambda^k(S)\ne V(H_{n,t})$, we must have
$v\notin\Lambda^k(S)$.
Thus, $\Lambda^k(S)=V(H_{n,t})\setminus\{v\}$ and
$|S\cap W|=d_S(v)\le k-1$.
Since $v\notin S$ and $|S\cap U|\le t$, we obtain
$|S|=|S\cap U|+|S\cap W|\le t+k-1$.
Together with $|S|\le n-1$, this yields
$|S|\le\min\{n-1,k+t-1\}$.
Therefore,
$$
\frac{|\Lambda^k(S)|}{|S|}
=\frac{n-1}{|S|}
\ge\frac{n-1}{\min\{n-1,k+t-1\}}
\ge\frac{n-1}{\min\{n,k+t\}}.
$$%
\end{case}%
This completes the proof.
\end{proof}

\Needspace{5\baselineskip}
We next determine the least integer $t$ for which
$\beta^k(H_{n,t})<r$. Recall that, for $(k,r)\ne(1,1)$ and
$n\ge(r+1)k$,
$$
p=p(n,k,r)=
\begin{cases}
1,&(r+1)k\leq n\le r(k+1)+1,\\[1mm]
\min\{n-(r+1)k+1,\ceil{\frac{n}{r}}-k\},&n\ge r(k+1)+2.
\end{cases}
$$

\begin{lem}\label{lem:p}
Let $k\ge1$, $r\ge1$ and $n\ge(r+1)k$ be integers with
$(k,r)\ne(1,1)$. Then $1\le p\le n-2$ and
$$
\beta^k(H_{n,p})<r.
$$
Moreover, $\beta^k(H_{n,t})\ge r$ for every integer $t$ with $1\le t<p$.
\end{lem}
\begin{proof}
Since $(k,r)\ne(1,1)$, we have $(r+1)k\ge3.$ Hence $n\ge3$.
Suppose first that $n\le r(k+1)+1$. Then $p=1\le n-2$.
By Lemma~\ref{lem:starbeta},
$\beta^k(H_{n,1})\le\frac{n-2}{k+1}\le r-\frac{1}{k+1}<r$.
Now suppose that $n\ge r(k+1)+2$.
Put $A=n-(r+1)k+1$ and $B=\ceil{\frac{n}{r}}-k$.
By the definition of $p$, we have $p=\min\{A,B\}$.
Since $n\ge(r+1)k$ and $n>r(k+1)$, we have $A\ge1$ and
$B\ge2$. Combining this with $(r+1)k\ge3$, we have
$1\le p\le A=n-(r+1)k+1\le n-2$.
For every integer $t$ with $1\le t\le n-2$, we have
$$
\begin{aligned}
&\frac{\max\{0,n-k-t\}}{k}<r
\iff n-k-t<rk
\iff t\ge A,\\
&\frac{n-1}{\min\{n,k+t\}}<r
\iff n\le r\min\{n,k+t\}
\iff n\le r(k+t)
\iff t\ge B.
\end{aligned}
$$
Since $p=\min\{A,B\}$, we have
$$
\min\left\{\frac{\max\{0,n-k-p\}}{k},
\frac{n-1}{\min\{n,k+p\}}\right\}<r.
$$
By Lemma~\ref{lem:starbeta}, we obtain that $\beta^k(H_{n,p})<r$.

Next we prove that $\beta^k(H_{n,t})\ge r$ for every integer $t$
with $1\le t<p$. Suppose to the contrary that
$\beta^k(H_{n,t})<r$ for some integer $t$ with $1\le t<p$.
Since $t<p=\min\{A,B\}$, we have $\frac{\max\{0,n-k-t\}}{k}\ge r$ and
$\frac{n-1}{\min\{n,k+t\}}\ge r.$
By Lemma~\ref{lem:starbeta}, we obtain
\begin{equation}\label{2.2.1}
\frac{n-1-\min\{t,k\}}{k+1}<r.
\end{equation}
Let $c=n-r(k+1)$. By $n\ge r(k+1)+2$, we have $c\geq2.$ Using (\ref{2.2.1}), we obtain that $\min\{t,k\}>c-1$.
Since $t$, $k$ and $c$ are integers, we have $2\le c\le\min\{t,k\}.$
If $r\ge2$, then $\frac{c}{r}\le\frac{c}{2}\le c-1$, and hence
$
B=1+\ceil{\frac{c}{r}}\le c.
$
If $r=1$, then $k\ge2$ and
$
A=c-k+2\le c.
$
In either case, we obtain
$$
p=\min\{A,B\}\le c\le\min\{t,k\}\le t<p,
$$
a contradiction. Hence $\beta^k(H_{n,t})\ge r$ for every integer $t$
with $1\le t<p$.
\end{proof}

For the next Lemmas \ref{edge:count} and \ref{edge:equality}, let $n,k,r$ be positive integers with $n\ge k$.
Let $S\subseteq V(G)$
be chosen as in the definition of $\beta^k(G)$ so that
$\frac{|\Lambda^k(S)|}{|S|}<r$. Define $T=V(G)\setminus\Lambda^k(S)$,
and let $|S|=s^*$ and $|T|=t^*$. We next establish lower bounds of $e(\overline G).$

\begin{lem}\label{edge:count} Let $G$ be a graph of order $n$ with $\beta^k(G)<r.$ Then\\
(i) $t^*\ge n-rs^*+1$.\\
(ii) If $s^*=k$, then $S\subseteq T$ and
$e(\overline G)\ge t^*-k\ge n-(r+1)k+1.$\\
(iii) If $s^*=k+x$ for some integer $x\ge1$, let $d=|T\setminus S|$.
Then
\begin{equation*}\label{eq:holes}
e(\overline G)\ge\max\left\{
x,\ceil{\frac{t^*x+d}{2}},t^*x-\binom{t^*}{2},
\ceil{\frac{(t^*-d)x}{2}}+d(x+1)
\right\}.
\end{equation*}
\end{lem}
\begin{proof}
(i) By $n-t^*=|\Lambda^k(S)|<rs^*$, we have $t^*>n-rs^*$.
Note that $n,t^*,r,s^*$ are integers. Then $t^*\ge n-rs^*+1$.

\textup{(ii)} Note that $s^*=k$. Then every vertex of $S$ has at most
$k-1$ neighbors in $S$, and hence $S\subseteq T$.
Let $v\in T\setminus S$. Suppose to the contrary that
$|\Lambda_G(v)\cap S|\ge k$. By the definition of $\Lambda^k(S)$,
we have $v\in\Lambda^k(S)$, which contradicts
$v\in T=V(G)\setminus\Lambda^k(S)$.
Hence $|\Lambda_G(v)\cap S|\le k-1$.
Since $v\notin S$ and $|S|=k$, we obtain
$
|\Lambda_{\overline G}(v)\cap S|
=|S|-|\Lambda_G(v)\cap S|
\ge k-(k-1)=1.
$
Since each edge of $\overline G$ joining $T\setminus S$ to $S$
has exactly one endpoint in $T\setminus S$, we obtain
$
e(\overline G)\ge|T\setminus S|=t^*-k\ge n-(r+1)k+1.
$

\textup{(iii)} Suppose that $s^*=k+x$, where $x\ge1$ is an integer.
Let $d=|T\setminus S|$. Then $|T\cap S|=t^*-d$.
For every $v\in T$, we have $|\Lambda_G(v)\cap S|\le k-1$.
Then
$
|\Lambda_{\overline G}(v)\cap S|\ge
\begin{cases}
x,&v\in T\cap S,\\
x+1,&v\in T\setminus S.
\end{cases}
$
Since $T\ne\emptyset$, we have $$e(\overline G)\ge x.$$
Summing over the vertices of $T$, we obtain
\begin{equation}\label{eq:sumlower}
t^*x+d=(t^*-d)x+d(x+1)
\le\sum_{v\in T}|\Lambda_{\overline G}(v)\cap S|.
\end{equation}
Each edge of $\overline G$ is counted at most twice, and only
edges of $\overline G[T\cap S]$ are counted twice.
Since $|T\cap S|\le t^*$, we obtain
$$
\begin{aligned}
\sum_{v\in T}|\Lambda_{\overline G}(v)\cap S|
\le\min\left\{2e(\overline G),e(\overline G)+e(\overline G[T\cap S])\right\}
\le\min\left\{2e(\overline G),e(\overline G)+\binom{t^*}{2}\right\}.
\end{aligned}
$$
Combining this with (\ref{eq:sumlower}), we obtain that
$$e(\overline G)\geq \ceil{\frac{t^*x+d}{2}}, e(\overline G)\ge t^*x+d-\binom{t^*}{2}\ge t^*x-\binom{t^*}{2}.$$

Recall that $|T\cap S|=t^*-d$ and each vertex of $T\cap S$ has
at least $x$ neighbors in $\overline G[S]$.
Since $T\cap S\subseteq S$,
$
2e(\overline G[S])
=\sum_{v\in S}|\Lambda_{\overline G}(v)\cap S|
\ge\sum_{v\in T\cap S}|\Lambda_{\overline G}(v)\cap S|
\ge(t^*-d)x.
$
Hence $e(\overline G[S])\ge\ceil{\frac{(t^*-d)x}{2}}$.
Moreover, there are at least $d(x+1)$ edges between $T\setminus S$
and $S$ in $\overline G$. These edges are distinct from the edges
of $\overline G[S]$. Hence
$$
e(\overline G)\ge e(\overline G[S])+d(x+1)
\ge\ceil{\frac{(t^*-d)x}{2}}+d(x+1).
$$
This completes the proof.
\end{proof}

\Needspace{5\baselineskip}
The following lemma describes the structure of $G$ in two equality
cases of Lemma~\ref{edge:count} (iii).

\begin{lem}\label{edge:equality}
Let $G$ be a graph of order $n$ with $\beta^k(G)<r,$ and let $s^*=k+x$ for an integer $x\ge1$.\\
(i) If $e(\overline G)=x$, then $G\cong K_n-K_{1,e(\overline G)}$.\\
(ii) If $2e(\overline G)=t^*x$, then $T\subseteq S$,
$e(\overline G[T])=e(\overline G)$, and $\overline G[T]$
is $x$-regular.
\end{lem}
\begin{proof}
By the definition of $T$, we have $|\Lambda_G(v)\cap S|\le k-1$
for every $v\in T$.

\textup{(i)} Note that $e(\overline G)=x$. Suppose to the contrary that
there exists a vertex $v\in T\setminus S$. Then
$$
e(\overline G)\ge|\Lambda_{\overline G}(v)\cap S|
=|S|-|\Lambda_G(v)\cap S|
\ge|S|-(k-1)=x+1=e(\overline G)+1,
$$
a contradiction. Hence $T\setminus S=\emptyset$, so $T\subseteq S$.
Since $T\ne\emptyset$, there exists a vertex $v\in T\subseteq S$.
Then
$$
e(\overline G)\ge|\Lambda_{\overline G}(v)\cap S|
=|S|-1-|\Lambda_G(v)\cap S|
\ge|S|-k=x=e(\overline G).
$$
Hence $|\Lambda_{\overline G}(v)\cap S|=e(\overline G)$.
It follows that all edges of $\overline G$ join $v$ to vertices of $S$,
and thus, $G\cong K_n-K_{1,e(\overline G)}$.

\textup{(ii)} Note that $2e(\overline G)=t^*x$.
By Lemma~\ref{edge:count} (iii), we have
$
t^*x+|T\setminus S|\le2e(\overline G)=t^*x,
$
so $T\setminus S=\emptyset$ and hence $T\subseteq S$.
For every $v\in T$, we have $|\Lambda_{\overline G}(v)\cap S|\ge x$.
Since each edge of $\overline G$ is counted at most twice, we obtain
$
2e(\overline G)=t^*x\le\sum_{v\in T}|\Lambda_{\overline G}(v)\cap S|\le2e(\overline G).
$
Hence
$
e(\overline G[T])=e(\overline G).
$
For every $v\in T$, we have
$$
|\Lambda_{\overline G[T]}(v)|
=|\Lambda_{\overline G}(v)\cap S|=x.
$$
Thus, $\overline G[T]$ is $x$-regular.
\end{proof}

\Needspace{5\baselineskip}
The following lemma describes the possible structures of graphs $G$
of order $n$ satisfying
$\beta^k(G)<r\le\beta^k(H_{n,e(\overline G)})$.

\begin{lem}\label{lem:classification}
Let $n,k,r$ be three positive integers with $n\ge(r+1)k$, and let
$G$ be a graph of order $n$ with $1\le e(\overline G)\le n-2$.
Suppose that $\beta^k(G)<r\le\beta^k(H_{n,e(\overline G)})$.
Then there exist a graph $F$, a vertex set $S\subseteq V(G)$,
and positive integers $x$ and $t^*$ such that $G\cong K_n-F$,
$|S|=k+x$, $|\Lambda^k(S)|=n-t^*<r(k+x)$,
and one of the following holds.\\
(i) $F\cong2K_2$, $k\in\{3,4\}$, $x=1$ and $t^*=4$.\\
(ii) $F\cong3K_2$, $k=5$, $x=1$ and $t^*=6$.\\
(iii) $F\cong K_3$, $k\in\{1,2,3\}$, $x=2$ and $t^*=3$.\\
(iv) $F\cong C_4$, $k=x=2$ and $t^*=4$.
\end{lem}
\begin{proof}
Since $\beta^k(G)<r$, there exists a vertex set $S\subseteq V(G)$
with $|S|\ge k$ such that $\Lambda^k(S)\ne V(G)$ and
$|\Lambda^k(S)|<r|S|$.
Define $T=V(G)\setminus\Lambda^k(S)$, and let $|S|=s^*$ and $|T|=t^*$.
Then $t^*\ge1$ and
\begin{equation}\label{eq:classification-ratio}
\frac{n-t^*}{s^*}<r.
\end{equation}
By $\beta^k(H_{n,e(\overline G)})\ge r\ge1$ and Lemma~\ref{lem:starbeta}, we have
$
r\le\beta^k(H_{n,e(\overline G)})
\le\frac{\max\{0,n-k-e(\overline G)\}}{k}.
$
Obviously, $\max\{0,n-k-e(\overline G)\}\ge rk>0$, so
$n-k-e(\overline G)>0$.
It follows that
\begin{equation}\label{eq:classification-starbound}
\frac{n-k-e(\overline G)}{k}
=\frac{\max\{0,n-k-e(\overline G)\}}{k}\ge r.
\end{equation}
We first show that $s^*\ge k+1$.
Suppose to the contrary that $s^*=k$.
By Lemma~\ref{edge:count} \textup{(ii)}, we have $e(\overline G)\geq t^*-k$.
Combining this with (\ref{eq:classification-starbound}), we obtain
$
r>\frac{n-t^*}{s^*}\ge\frac{n-k-e(\overline G)}{k}\ge r,
$
a contradiction. Hence $|S|=s^*=k+x$ for some integer $x\ge1$. By the choice of $S$ and the definition of $T$, we have
$|\Lambda^k(S)|=n-t^*<r|S|=r(k+x)$.

By Lemma~\ref{edge:count} \textup{(iii)}, we have
$e(\overline G)\ge x$ and $e(\overline G)\ge\ceil{\frac{t^*x}{2}}$.
If $e(\overline G)=x$, then Lemma~\ref{edge:equality} \textup{(i)} implies
$G\cong H_{n,e(\overline G)}$, and so $\beta^k(G)=\beta^k(H_{n,e(\overline G)})\ge r$,
a contradiction. Thus, $e(\overline G)>x$.
By (\ref{eq:classification-ratio}), $s^*=k+x$ and Lemma~\ref{lem:starbeta}, we obtain
$$
\frac{n-t^*}{k+x}<r\le
\min\left\{\frac{n-k-e(\overline G)}{k},\frac{n-1}{k+e(\overline G)}\right\}.
$$
It follows that
$xn>(k+x)(k+e(\overline G))-kt^*$ and $(e(\overline G)-x)n<t^*(k+e(\overline G))-(k+x).$
Since $x>0$ and $e(\overline G)-x>0$, we have
$
(e(\overline G)-x)\bigl((k+x)(k+e(\overline G))-kt^*\bigr)
<x\bigl(t^*(k+e(\overline G))-(k+x)\bigr).
$
Combining this with $kt^*(e(\overline G)-x)+xt^*(k+e(\overline G))=t^*\,e(\overline G)(k+x)$, we obtain
$(k+x)((e(\overline G)-x)(k+e(\overline G))+x)<t^*\,e(\overline G)(k+x)$.
Hence
$
t^*\,e(\overline G)>(e(\overline G)-x)(k+e(\overline G))+x=e(\overline G)+(e(\overline G)-x)(k+e(\overline G)-1),
$
which implies
\begin{equation}\label{eq:key}
t^*>1+\frac{(e(\overline G)-x)(k+e(\overline G)-1)}{e(\overline G)}.
\end{equation}
Since $k\ge1$, we have $\frac{k+e(\overline G)-1}{e(\overline G)}\ge1$.
Combining this with (\ref{eq:key}), we have $t^*>1+e(\overline G)-x$.
As $t^*,e(\overline G),x$ are integers, it follows that
\begin{equation}\label{eq:integerbound}
e(\overline G)\le t^*+x-2.
\end{equation}

\begin{claim}\label{clm:classification-size}
$t^*\le s^*$.
\end{claim}
\begin{proof}
Suppose to the contrary that $t^*=s^*+d$ for some integer $d\ge1$.
Let $d'=|T\setminus S|$. Then
$d'=t^*-|T\cap S|\ge t^*-s^*=d$.
By Lemma~\ref{edge:count} \textup{(iii)}, $d'\ge d$ and $t^*=s^*+d,$ we have
\begin{equation}\label{eq:outside}
\begin{aligned}
e(\overline G)\ge\ceil{\frac{(t^*-d')x}{2}}+d'(x+1)
 =\ceil{\frac{t^*x+d'(x+2)}{2}}
 \ge\ceil{\frac{s^*x}{2}}+d(x+1).
\end{aligned}
\end{equation}
According to (\ref{eq:integerbound}), we have
$e(\overline G)\le t^*+x-2=s^*+d+x-2$.
Combining this with (\ref{eq:outside}), we obtain
$s^*x+2d(x+1)\le2s^*+2d+2x-4$, and hence
\begin{equation}\label{eq:claim-sizebound}
(x-2)(s^*-2)+2dx\le0.
\end{equation}

If $x\ge2$, then $s^*=k+x\ge3$ and
$(x-2)(s^*-2)+2dx>0$. This contradicts (\ref{eq:claim-sizebound}).
Thus, $x=1$ and $s^*=k+1\ge2$.
If $s^*=2$, by (\ref{eq:integerbound}) and (\ref{eq:outside}), we obtain that
$1+2d\le e(\overline G)\le d+1$, contrary to $d\ge1$.
Hence $s^*\ge3$.
By (\ref{eq:outside}), we have
$e(\overline G)\ge\ceil{\frac{s^*}{2}}+2d\ge2d+2\ge4$.
Substituting $k=s^*-1$ and $t^*=s^*+d$ into (\ref{eq:key}), we obtain
$e(\overline G)(s^*+d-1)>(e(\overline G)-1)(s^*+e(\overline G)-2)$, and hence
$s^*>e(\overline G)^2-e(\overline G)(d+2)+2$.
By (\ref{eq:outside}), we also have $s^*\le2e(\overline G)-4d$.
Combining $s^*>e(\overline G)^2-e(\overline G)(d+2)+2$,
$s^*\le2e(\overline G)-4d$ and $e(\overline G)\ge2d+2\ge4$, we obtain
$
0>e(\overline G)^2-e(\overline G)(d+4)+4d+2=(e(\overline G)-4)(e(\overline G)-d)+2\ge2,
$
a contradiction.
Therefore, $t^*\le s^*=k+x$.
\end{proof}
\begin{claim}\label{clm:classification-equality}
$2e(\overline G)=t^*x$.
\end{claim}
\begin{proof}
We consider the cases $t^*\le x+1$ and $t^*\ge x+2$.

\resetcases
\begin{case}{$t^*\le x+1$.}\label{case:classification-small}

By $e(\overline G)>x$ and (\ref{eq:integerbound}), we have
$x<e(\overline G)\le t^*+x-2$. This implies that $t^*\ge3$.
Combining Lemma~\ref{edge:count} \textup{(iii)} with
(\ref{eq:integerbound}), we obtain
\begin{equation}\label{eq:classification-case1bounds}
t^*x-\binom{t^*}{2}\le e(\overline G)\le t^*+x-2.
\end{equation}
Since $x\ge t^*-1$, we have
$
(t^*-1)^2\le(t^*-1)x\le\binom{t^*}{2}+t^*-2,
$
which implies $(t^*-2)(t^*-3)\le0$.
According to $t^*\ge3$, we know that $t^*=3$.
Substituting $t^*=3$ into (\ref{eq:classification-case1bounds}), we obtain
$3x-3\le e(\overline G)\le x+1$.
Since $x\ge t^*-1=2$, we obtain $x=2$ and $2<e(\overline G)\le3$.
Hence $e(\overline G)=3$.
By (\ref{eq:key}), we have $3>1+\frac{k+2}{3}$,
and so $k\in\{1,2,3\}$.%
\end{case}

\begin{case}{$t^*\ge x+2$.}\label{case:classification-large}

By Lemma~\ref{edge:count} \textup{(iii)} and
(\ref{eq:integerbound}), we have
$\frac{t^*x}{2}\le e(\overline G)\le t^*+x-2$, which implies
$(t^*-2)(x-2)\le0$.
Since $t^*\ge3$ and $x\ge1$, we have $x=1$ or $x=2$.
We first consider $x=2$.
Then $t^*\ge4$ and $t^*\le e(\overline G)\le t^*.$ This implies that $e(\overline G)=t^*$.
According to $t^*\le s^*=k+2$, we have $k\ge t^*-2$.
By (\ref{eq:key}), we obtain
$
t^*(t^*-1)>(t^*-2)(k+t^*-1)\ge(t^*-2)(2t^*-3).
$
This implies that $(t^*)^2-6t^*+6<0$.
Note that $t^*$ is an integer with $t^*\ge4$. Then $t^*=4$.
Substituting $t^*=e(\overline G)=4$ and $x=2$ into (\ref{eq:key}), we have
$4>1+\frac{k+3}{2}$.
Combining $k\ge t^*-2=2$, this implies $k=2$.

Next we consider $x=1$.
According to $3\le t^*\le s^*=k+1$, we have $k\ge t^*-1$.
By (\ref{eq:key}), we obtain
$$
t^*>1+\frac{(e(\overline G)-1)(k+e(\overline G)-1)}{e(\overline G)}
\ge1+\frac{(e(\overline G)-1)(t^*+e(\overline G)-2)}{e(\overline G)},
$$
and hence $t^*>(e(\overline G)-1)^2+1$.
By Lemma~\ref{edge:count} \textup{(iii)}, we have $t^*\le2e(\overline G)$.
Thus, $(e(\overline G)-2)^2<2$.
Since $e(\overline G)>x=1$ and $e(\overline G)$ is an integer, we obtain $e(\overline G)=2$ or $3$.

If $e(\overline G)=2$, then $3\le t^*\le4$.
Suppose that $t^*=3$.
By (\ref{eq:key}), we have $3>1+\frac{k+1}{2}$.
Combining $k\ge t^*-1=2$, we have $k=2$.
By Lemma~\ref{lem:starbeta}, we obtain
$$
r\le\beta^k(H_{n,e(\overline G)})\le\frac{n-3}{3}
=\frac{n-t^*}{k+x}<r,
$$
a contradiction. Hence $t^*=4,$ and $k\ge t^*-1=3$.
By (\ref{eq:key}) and $4>1+\frac{k+1}{2}$, we have $k=3$ or $k=4$.

If $e(\overline G)=3$, then $5<t^*\le6$. This implies that $t^*=6$.
Since $t^*\le k+1$, we have $k\ge5$.
By (\ref{eq:key}), we obtain $6>1+\frac{2(k+2)}{3}$,
which implies $k\le5$.
Thus, $k=5$.%
\end{case}
By Cases \ref{case:classification-small} and \ref{case:classification-large}, we have
$
(x,t^*,e(\overline G))
\in\{(2,3,3),(2,4,4),(1,4,2),(1,6,3)\}.
$
Thus, $2e(\overline G)=t^*x$.
\end{proof}
By Claim~\ref{clm:classification-equality} and
Lemma~\ref{edge:equality} \textup{(ii)}, we obtain
$
T\subseteq S,$ $e(\overline G[T])=e(\overline G),
$
and $\overline G[T]$ is $x$-regular.
Taking $F=\overline G[T]$, we have $G\cong K_n-F$.
If $x=1$, then every component of $F$ is $K_2$.
Hence $F\cong2K_2$ for $t^*=4$ and $F\cong3K_2$ for $t^*=6$.
If $x=2$, then every component of $F$ is a cycle of length at least three.
Thus, $F\cong K_3$ for $t^*=3$ and $F\cong C_4$ for $t^*=4$.
This proves \textup{(i)-(iv)}. This completes the proof.
\end{proof}

The following lemma shows that, when $(k,r)\ne(1,1)$,
the possibilities in Lemma~\ref{lem:classification} reduce
to a single exceptional graph.

\begin{lem}\label{lem:exception}
Let $n, k, r$ be positive integers with $n\ge(r+1)k$ and
$(k,r)\ne(1,1)$, and let $G$ be a graph of order $n$
with $1\le e(\overline G)\le n-2$.
If $\beta^k(G)<r\le\beta^k(H_{n,e(\overline G)}),$
then $(n,k,r)=(11,3,2)$ and $G\cong K_{11}-2K_2$.
\end{lem}
\begin{proof}
By Lemma~\ref{lem:classification}, we have $G\cong K_n-F$ and
$|\Lambda^k(S)|=n-t^*<r(k+x)$, where $|S|=k+x$.
Moreover, $F,k,x,t^*$ satisfy one of
\textup{(i)-(iv)} in Lemma~\ref{lem:classification}.
By $\beta^k(H_{n,e(\overline G)})\ge r\ge1$ and
Lemma~\ref{lem:starbeta}, we have
$
\frac{\max\{0,n-k-e(\overline G)\}}{k}\ge r.
$
Hence $\max\{0,n-k-e(\overline G)\}\ge rk>0$. This implies that
$n-k-e(\overline G)\ge rk>0$.
Thus, $k+e(\overline G)<n$.
Applying Lemma~\ref{lem:starbeta}, we obtain
$$
r\le\min\left\{
\frac{n-k-e(\overline G)}{k},
\frac{n-1}{k+e(\overline G)},
\frac{n-1-\min\{e(\overline G),k\}}{k+1}
\right\}.
$$
By the three lower bounds above, we obtain
$n\geq \max\left\{(r+1)k+e(\overline G),\,
r(k+e(\overline G))+1,\,
r(k+1)\right.$
$\left.+1+\min\{e(\overline G),k\}\right\}=\theta.$ Since $n-t^*<r(k+x)$,
we have $n\le r(k+x)+t^*-1$. This implies that
\begin{equation}\label{eq:integer_n}
\begin{aligned}
&\theta
\le n\le r(k+x)+t^*-1.
\end{aligned}
\end{equation}

\resetcases
\begin{case}{$F\cong2K_2$.}\label{case:exception-twoedges}

By Lemma~\ref{lem:classification} \textup{(i)}, we have
$k\in\{3,4\}$, $x=1$ and $t^*=4$.
By $F\cong2K_2$, we have $e(\overline G)=e(F)=2$.
We first consider $k=4$.
By (\ref{eq:integer_n}), we have
$4r+6\le5r+3$ and $6r+1\le5r+3$,
which imply $3\le r\le2$, a contradiction.
Hence $k=3$.
Substituting these values into (\ref{eq:integer_n}), we obtain
$
\max\{3r+5,5r+1,4r+3\}\le n\le4r+3.
$
In particular, $3r+5\le4r+3$ and $5r+1\le4r+3$,
which imply $2\le r\le2$.
Thus, $r=2$.
Substituting $r=2$ into the bounds on $n$, we have $11\le n\le11$,
and hence $n=11$.%
\end{case}

\begin{case}{$F\cong3K_2$.}\label{case:exception-threeedges}

By Lemma~\ref{lem:classification} \textup{(ii)}, we have
$k=5$, $x=1$ and $t^*=6$.
By $F\cong3K_2$, we have $e(\overline G)=e(F)=3$.
Substituting these values into (\ref{eq:integer_n}), we obtain
$5r+8\le6r+5$ and $8r+1\le6r+5$.
It follows that $3\le r\le2$, a contradiction.%
\end{case}

\begin{case}{$F\cong K_3$.}\label{case:exception-triangle}

By Lemma~\ref{lem:classification} \textup{(iii)}, we have
$k\in\{1,2,3\}$, $x=2$ and $t^*=3$.
By $F\cong K_3$, we have $e(\overline G)=e(F)=3$.
If $k=1$, then (\ref{eq:integer_n}) gives $4r+1\le3r+2$.
Hence $r\le1$.
Together with $r\ge1$, this implies $r=1$,
contrary to $(k,r)\ne(1,1)$.
If $k=2$, then (\ref{eq:integer_n}) gives
$2r+5\le4r+2$ and $5r+1\le4r+2$.
Since $r$ is an integer, we obtain $2\le r\le1$, a contradiction.
If $k=3$, then (\ref{eq:integer_n}) gives
$3r+6\le5r+2$ and $6r+1\le5r+2$,
which imply $2\le r\le1$, a contradiction.%
\end{case}

\begin{case}{$F\cong C_4$.}\label{case:exception-cycle}

By Lemma~\ref{lem:classification} \textup{(iv)}, we have
$k=x=2$ and $t^*=4$.
By $F\cong C_4$, we have $e(\overline G)=e(F)=4$.
Substituting these values into (\ref{eq:integer_n}), we obtain
$2r+6\le4r+3$ and $6r+1\le4r+3$.
Since $r$ is an integer, we have $2\le r\le1$, a contradiction.%
\end{case}

Cases \ref{case:exception-threeedges}-\ref{case:exception-cycle} all lead to contradictions, and Case \ref{case:exception-twoedges} with $k=4$
is also impossible.
Thus, only Case \ref{case:exception-twoedges} with $k=3$ remains, in which $r=2$ and $n=11$.
Therefore, $(n,k,r)=(11,3,2)$ and $F\cong2K_2$.
Hence $G\cong K_{11}-2K_2$.
\end{proof}

The following property of binding sets will be used in the proofs
of Theorems~\ref{edge:first} and \ref{thm:first}.

\begin{lem}[Kane, Mohanty and Straus \cite{Kane1981}]\label{lem:bindingset}
Let $G$ be a graph with $b(G)<1$.
Then every binding set of $G$ is independent.
\end{lem}

For integers $n\ge3$ and $1\le t\le\floor{\frac{n-1}{2}}$, let
$
F_{n,t}=K_t\vee\bigl(K_{n-2t-1}\cup(t+1)K_1\bigr),
$
where $K_0$ denotes the graph with no vertices.
We now show that the maximum size and spectral radius in
$\mathcal G_{n,1,1}$ are attained within this family.

\begin{lem}\label{lem:family}
Let $G$ be a graph of order $n$ without isolated vertices
and with $b(G)<1$. Then\\
(i) There exists an integer $t$ with
$1\le t\le\floor{\frac{n-1}{2}}$ such that a graph isomorphic
to $F_{n,t}$ can be obtained from $G$ by adding edges.
In particular, $n\ge3$.\\
(ii) $F_{n,t}\in\mathcal G_{n,1,1}$ for every integer $t$
with $1\le t\le\floor{\frac{n-1}{2}}$.
\end{lem}
\begin{proof}
\textup{(i)} Let $S$ be a binding set of $G$.
Put $|S|=s^*$ and $|\Lambda_G(S)|=\Lambda$.
Since $G$ has no isolated vertices and $b(G)=\frac{\Lambda}{s^*}<1$,
we have $1\le\Lambda\le s^*-1$.
By Lemma~\ref{lem:bindingset}, $S$ is independent, and hence
$S\cap\Lambda_G(S)=\emptyset$.
Thus, $n\ge s^*+\Lambda\ge2\Lambda+1\ge3$. This implies that
$1\le\Lambda\le\floor{\frac{n-1}{2}}$.
Choose a set $S_1\subseteq S$ with $|S_1|=\Lambda+1$, and put
$S_2=V(G)\setminus(\Lambda_G(S)\cup S_1)$.
Then $\Lambda_G(S),S_1,S_2$ are pairwise disjoint, with union $V(G)$ and sizes
$\Lambda,\Lambda+1,n-2\Lambda-1$, respectively.
Since $S$ is independent and $S_1\subseteq S$, we have
$E(G[S_1])=\emptyset$ and $\Lambda_G(S_1)\subseteq\Lambda_G(S)$.
Hence there are no edges between $S_1$ and $S_2$.
Let $G'$ be the graph obtained from $G$ by adding all missing edges
within $\Lambda_G(S)\cup S_2$ and between $\Lambda_G(S)$ and $S_1$.
Then $G'[\Lambda_G(S)\cup S_2]$ is complete, $S_1$ is independent in $G'$,
and every vertex of $S_1$ is adjacent precisely to the vertices of $\Lambda_G(S)$.
Thus, $G'\cong K_\Lambda\vee\bigl(K_{n-2\Lambda-1}\cup(\Lambda+1)K_1\bigr)=F_{n,\Lambda}$.
Taking $t=\Lambda$, we obtain Lemma~\ref{lem:family} \textup{(i)}.

\vspace{1mm}
\noindent\textup{(ii)} Let $t$ be an integer with
$1\le t\le\floor{\frac{n-1}{2}}$.
Let $S_1$ be the vertex set of the $(t+1)K_1$ part of $F_{n,t}$.
By the definition of $F_{n,t}$, we have that
$\Lambda_{F_{n,t}}(S_1)$ is the vertex set of the $K_t$ part.
Thus, $\Lambda_{F_{n,t}}(S_1)\ne V(F_{n,t})$ and
$b(F_{n,t})\le\frac{|\Lambda_{F_{n,t}}(S_1)|}{|S_1|}
=\frac{t}{t+1}<1$.
Combining this with $\delta(F_{n,t})=t\ge1$, we obtain
$F_{n,t}\in\mathcal G_{n,1,1}$.
This proves Lemma~\ref{lem:family} \textup{(ii)}.
\end{proof}

\section{Proofs of Theorems~\ref{edge:main} and \ref{edge:first}}\label{edge:proofs}

\begin{proof}[\normalfont\bfseries Proof of Theorem~\ref{edge:main}]
Let $G^*$ be a graph with maximum size in $\mathcal G_{n,k,r}$.
Since $G^*\in\mathcal G_{n,k,r}$, we have
$\beta^k(G^*)<r\le\beta^k(K_n)$, and hence $e(\overline{G^*})\ge1$.
Choose a vertex set $S\subseteq V(G^*)$ such that
$|S|\ge k$, $\Lambda_{G^*}^k(S)\ne V(G^*)$ and
$\frac{|\Lambda_{G^*}^k(S)|}{|S|}<r$.
Let $T=V(G^*)\setminus\Lambda_{G^*}^k(S)$,
$|S|=s^*=k+x$, $|T|=t^*$ and $d=|T\setminus S|$.
Then $x\ge0$ and $t^*\ge1$.
By Lemma~\ref{edge:count} \textup{(i)}, we have
\begin{equation}\label{edge:equalwitness}
t^*\ge n-r(k+x)+1.
\end{equation}
By Lemma~\ref{lem:p}, we have $1\le p\le n-2$ and
$\beta^k(H_{n,p})<r$.
Since $\delta(H_{n,p})=n-p-1\ge1$, we obtain
$H_{n,p}\in\mathcal G_{n,k,r}$, and hence
$
e(G^*)\ge e(H_{n,p})=\binom n2-p.
$
Thus, $1\le e(\overline{G^*})\le p$.

\vspace{1mm}
\noindent\textup{(i)} Note that $(n,k,r)=(11,3,2)$.
Then we claim that $e(\overline{G^*})<p$.
Suppose to the contrary that $e(\overline{G^*})\ge p$.
By $1\leq e(\overline{G^*})\le p$, we have
$e(G^*)=\binom{11}{2}-p=\binom{11}{2}-3=52$.
Let $G'=K_{11}-2K_2$, and let $S_0$ be the set of four vertices incident with the deleted edges. Then
$$
\beta^3(G')\le\frac{|\Lambda_{G'}^3(S_0)|}{|S_0|}=\frac74<2.
$$
By $\delta(G')\geq1$, we have $G'\in\mathcal G_{11,3,2}$.
However,
$
e(G^*)=52<53=e(G'),
$
which contradicts the maximality of $G^*$.
Thus, $e(\overline{G^*})<p$.
According to Lemma~\ref{lem:p}, we have
$
\beta^k(G^*)<r\le\beta^k(H_{n,e(\overline{G^*})}).
$
By Lemma~\ref{lem:exception}, we obtain $G^*\cong K_{11}-2K_2$.
Thus, $e(G)\le e(G^*)=53$ with equality if and only if
$G\cong K_{11}-2K_2$.
This proves Theorem~\ref{edge:main} \textup{(i)}.

\vspace{1mm}
\noindent\textup{(ii)} Note that $(n,k,r)\ne(11,3,2)$.
If $e(\overline{G^*})<p$, by Lemmas~\ref{lem:p} and \ref{lem:exception},
we obtain that $(n,k,r)=(11,3,2)$, a contradiction. Hence $e(\overline{G^*})\geq p.$
Note that $e(\overline{G^*})\le p$. Then $e(\overline{G^*})=p$.
Let $F$ be the subgraph of $\overline{G^*}$ with
$E(F)=E(\overline{G^*})$ and no isolated vertices.
Thus, $G^*=K_n-F$ and $e(F)=p$.
The labels \textup{(a)-(h)} below refer to the conditions
in the definition of $\mathcal F_{n,k,r}$.
If $p=1$, then $F\cong K_{1,p}$ and \textup{(a)} holds.
Next we consider $p\ge2$.
Put $A=n-(r+1)k+1$ and $B=\ceil{\frac{n}{r}}-k$.
By the definition of $p$, we have $n\ge r(k+1)+2$ and
$p=\min\{A,B\}$.
Thus,
\begin{equation}\label{edge:equalbounds}
n\ge\max\{(r+1)k+p-1,\,r(k+p-1)+1\}.
\end{equation}
Recall that $|S|=s^*=k+x.$
If $x\ge1$, by Lemma~\ref{edge:count} \textup{(iii)}, we have
$p=e(\overline{G^*})\geq x$.
The inequality $x\le p$ also holds trivially when $x=0$.
Hence $0\le x\le p$.

\begin{claim}\label{clm:edge-support}
For every $uv\in E(\overline{G^*})$, we have
$
u\in T,\ v\in S\ \text{or}\ u\in S,\ v\in T.
$
\end{claim}
\begin{proof}
Suppose to the contrary that there exists an edge $uv\in E(\overline{G^*})$
which does not join $T$ to $S$.
Let $G'=G^*+uv$.
Then $|\Lambda_{G'}(w)\cap S|=|\Lambda_{G^*}(w)\cap S|$
for every $w\in T$, and hence
$T\subseteq V(G')\setminus\Lambda_{G'}^k(S)$.
Since $T\ne\emptyset$, we have $\Lambda_{G'}^k(S)\ne V(G')$ and
$
\beta^k(G')\le\frac{|\Lambda_{G'}^k(S)|}{|S|}
\le\frac{n-t^*}{k+x}<r.
$
Combining this with $\delta(G')\ge\delta(G^*)\ge1$, we have
$G'\in\mathcal G_{n,k,r}$.
Note that $e(G')=e(G^*)+1$. This contradicts the maximality of $G^*$.
\end{proof}
\noindent By Claim~\ref{clm:edge-support}, we obtain that
\begin{equation}\label{edge:sumexact}
\sum_{v\in T}|\Lambda_{\overline{G^*}}(v)\cap S|
=e(\overline{G^*})+e(\overline{G^*}[T\cap S])=p+e(\overline{G^*}[T\cap S]).
\end{equation}

If $x=p$, by Lemma~\ref{edge:equality} \textup{(i)}, we have
$F\cong K_{1,p}$ and \textup{(a)} holds.
Next we consider the three cases $x=0$, $x=1$ and $2\le x<p$.

\resetcases
\begin{case}{$x=0$.}\label{case:edge-zero}

Then $s^*=k$. By Lemma~\ref{edge:count} \textup{(ii)}, we have
$S\subseteq T$ and
$
e(\overline{G^*})=p\ge t^*-k\ge n-(r+1)k+1=A\ge p.
$
Thus, $p=A=t^*-k$.
Note that every vertex of $T\setminus S$ has at least one neighbor in $S$
in $\overline{G^*}$.
By $e(\overline{G^*})=|T\setminus S|=p$, we obtain that every vertex of
$T\setminus S$ has degree one in $\overline{G^*}$, and every edge
of $\overline{G^*}$ joins $T\setminus S$ to $S$.
Hence $F$ is a star forest with at most $|S|=k$ components.
According to $n=(r+1)k+p-1$, we have \textup{(b)}.%
\end{case}

\begin{case}{$x=1$.}\label{case:edge-one}

We first claim that $T\subseteq S$.
Suppose to the contrary that $d=|T\setminus S|\ge1$.
By Lemma~\ref{edge:count} \textup{(iii)}, we have
\begin{equation}\label{edge:xonecount}
2p=2e(\overline{G^*})\ge2\ceil{\frac{t^*-d}{2}}+4d\ge t^*+3d.
\end{equation}
We first consider $r\ge2$.
By (\ref{edge:equalwitness}) and (\ref{edge:equalbounds}),
we have $t^*\ge r(p-2)+2\ge2p-2$.
Combining this with (\ref{edge:xonecount}), we obtain that
$
2p\ge t^*+3d\ge2p+1,
$
a contradiction.
Next we consider $r=1$.
Since $(k,r)\ne(1,1)$, we have $k\ge2$.
By (\ref{edge:equalwitness}) and (\ref{edge:equalbounds}), we have
$t^*\ge k+p-1$.
Since $|S|=k+1$, $d\ge t^*-|S|\ge p-2$.
By (\ref{edge:xonecount}), we have
$
2p\ge k+p-1+3\max\{1,p-2\}.
$
If $p=2$, then $4\ge k+4\ge6$, a contradiction.
If $p\ge3$, then $6\le2p\le7-k\le5$, a contradiction.
Hence $d=0$. This implies that $T\subseteq S$.

By Claim~\ref{clm:edge-support} and $T\subseteq S$,
we have $e(\overline{G^*}[S])=p$.
According to $T\subseteq S$, $|S|=k+1$ and the definition of $T$, we have
$
T=\{v\in S:d_{\overline{G^*}[S]}(v)\ge1\}.
$
Thus, $V(F)=T$ and $t^*=|V(F)|\ge2$.

\begin{claim}\label{clm:edge-minimum}
$p=\ceil{\frac{t^*}{2}}$.
\end{claim}
\begin{proof}
Since $V(F)=T$ and $F$ has no isolated vertices, we have
$
2p=\sum_{v\in T}d_F(v)\ge t^*,
$
and hence $p\ge\ceil{\frac{t^*}{2}}$.
Suppose to the contrary that $p>\ceil{\frac{t^*}{2}}$.
Construct a graph $F'$ on $T$ such that
$$
F'\cong
\begin{cases}
\frac{t^*}{2}K_2,&t^*\text{ is even},\\[1mm]
P_3\cup\frac{t^*-3}{2}K_2,&t^*\text{ is odd}.
\end{cases}
$$
Then $F'$ has no isolated vertices and $e(F')=\ceil{\frac{t^*}{2}}$.
Let $G'=K_n-F'$.
Since $T\subseteq S$, $|S|=k+1$ and $F'$ has no isolated vertices,
we obtain that
$
\Lambda_{G'}^k(S)=V(G')\setminus T\ne V(G').
$
Hence
$
\beta^k(G')\le\frac{|\Lambda_{G'}^k(S)|}{|S|}
=\frac{n-t^*}{k+1}<r.
$
Moreover,
$
\delta(G')\ge n-1-e(F')
=n-1-\ceil{\frac{t^*}{2}}\ge n-1-p\ge1.
$
Thus, $G'\in\mathcal G_{n,k,r}$.
However,
$
e(G')=\binom n2-\ceil{\frac{t^*}{2}}>\binom n2-p=e(G^*),
$
which contradicts the maximality of $G^*$.
\end{proof}

By Claim~\ref{clm:edge-minimum}, we have $t^*=2p$ or $t^*=2p-1$.
Recall that $V(F)=T$ and $|S|=k+1$.
By $T\subseteq S$ and (\ref{edge:equalwitness}), we obtain
$t^*\le k+1$ and $n\le r(k+1)+t^*-1.$
We first consider $t^*=2p$.
Since $F$ has no isolated vertices and
$\sum_{v\in V(F)}d_F(v)=2p=|V(F)|$, every vertex of $F$
has degree one.
Hence $F\cong pK_2$.
By $t^*=2p$, we have
$2p\le k+1$ and $n\le r(k+1)+2p-1.$
Next we consider $t^*=2p-1$.
Since $F$ has no isolated vertices and
$\sum_{v\in V(F)}d_F(v)=2p=|V(F)|+1$, exactly one vertex of $F$
has degree two and all the others have degree one.
Thus, $F\cong P_3\cup(p-2)K_2$.
By $t^*=2p-1$, we obtain
$2p-1\le k+1,$ and $n\le r(k+1)+2p-2.$

For $p=2$, we have $F\cong2K_2$ or $P_3$, which satisfy
\textup{(c)} and \textup{(a)}, respectively.
For $p=3$, we have $F\cong3K_2$ or $P_3\cup K_2$, which satisfy
\textup{(e)} and \textup{(d)}, respectively.
Now suppose that $p\ge4$.
We show that $A=p$.
Suppose to the contrary that $A>p$.
Then $n\ge(r+1)k+p$.
Combining this with $t^*\le k+1$, $n\le r(k+1)+t^*-1$, we obtain
$
t^*-1\le k\le r+t^*-p-1,
$
and hence $r\ge p$.
By (\ref{edge:equalbounds}), we have
$
r(k+p-1)+1\le n\le r(k+1)+t^*-1.
$
This implies $r(p-2)\le t^*-2$.
Since $t^*\le2p$ and $p\ge4$, we obtain
$
2p-2\ge t^*-2\ge r(p-2)\ge p(p-2)>2p-2,
$
a contradiction.
Therefore, $A=p$, and hence $n=(r+1)k+p-1$.
The graphs $pK_2$ and $P_3\cup(p-2)K_2$ have $p$ and $p-1$
components, respectively.
By $k\ge2p-2\ge p$, we obtain that $F\in\mathcal S_{k,p}$.
Thus, \textup{(b)} holds.%
\end{case}

\begin{case}{$2\le x<p$.}\label{case:edge-large}

\begin{claim}\label{clm:edge-parameters}
$x=2$ and $(p,t^*)\in\{(3,2),(3,3),(4,4)\}$.
\end{claim}
\begin{proof}
Let $y=p-x\ge1$.
By Lemma~\ref{edge:count} \textup{(iii)}, we have
\begin{equation}\label{edge:parametercount}
t^*x\le2p=2(x+y),\qquad t^*x-\binom{t^*}{2}\le p.
\end{equation}
Combining (\ref{edge:equalwitness}) and (\ref{edge:equalbounds}),
we obtain
\begin{equation}\label{edge:parameterlower}
t^*\ge\max\{k+p-rx,\,r(y-1)+2\}.
\end{equation}

We first consider $r\ge2$.
By $y\ge1$ and (\ref{edge:parameterlower}), we have
$t^*\ge r(y-1)+2\ge2y$.
Using (\ref{edge:parametercount}), we obtain that
$
2xy\le t^*x\le2(x+y),
$
and hence $(x-1)(y-1)\le1$.
If $y\ge2$, then $x,y\ge2$ and $(x-1)(y-1)=1$.
Hence $x=y=2$ and $p=4$.
By (\ref{edge:parametercount}) and (\ref{edge:parameterlower}),
we have $4\le t^*\le4$.
Thus, $(p,t^*)=(4,4)$.
Next we assume that $y=1$.
Then $p=x+1$.
By (\ref{edge:parametercount}) and (\ref{edge:parameterlower}), we have
$
2\le t^*\le\frac{2p}{x}=2+\frac2x.
$
If $x\ge3$, then $t^*=2$.
Combining this with (\ref{edge:parametercount}), we obtain
$
x+1=p\ge t^*x-\binom{t^*}{2}=2x-1\ge x+2,
$
a contradiction.
Hence $x=2$, $p=3$ and $2\le t^*\le3$.
Thus, $(p,t^*)\in\{(3,2),(3,3)\}$.

We now consider $r=1$.
Since $(k,r)\ne(1,1)$, we have $k\ge2$.
By (\ref{edge:parameterlower}), we obtain
$t^*\ge k+p-x=k+y$.
Combining this with (\ref{edge:parametercount}), we have
$(k+y)x\le2(x+y)$, and hence
$
0\le(x-2)y\le x(2-k)\le0.
$
Since $y\ge1$ and $x\ge2$, it follows that $x=k=2$.
By (\ref{edge:parametercount}) and (\ref{edge:parameterlower}), we obtain
$
p=y+2\le t^*\le\frac{2p}{x}=p.
$
Thus, $t^*=p$ and $t^*x=2p=2e(\overline{G^*})$.
By Lemma~\ref{edge:equality} \textup{(ii)}, we have $T\subseteq S$.
Since $p>x=2$, we obtain
$
3\le p=t^*\le|S|=k+x=4.
$
Hence $(p,t^*)\in\{(3,3),(4,4)\}$.
\end{proof}

By Claim~\ref{clm:edge-parameters}, we have $x=2$ and
$(p,t^*)\in\{(3,2),(3,3),(4,4)\}$.
Recall that $|S|=k+2$ and
$|\Lambda_{G^*}^k(S)|=n-t^*<r(k+2)$.
Note that the parameters are integers. Then $n\le r(k+2)+t^*-1$.
We determine $F$ for each of the three parameter pairs.

\begin{subcase}{$(p,t^*)=(3,2)$.}\label{subcase:edge-path}

By (\ref{eq:sumlower}) and (\ref{edge:sumexact}), we have
$$
4+d\le\sum_{v\in T}|\Lambda_{\overline{G^*}}(v)\cap S|
=3+e(\overline{G^*}[T\cap S])\le4.
$$
Hence $d=0$, $T\subseteq S$ and $e(\overline{G^*}[T])=1$.
Let $T=\{u,v\}$. Then $uv\in E(\overline{G^*})$.
Each of $u$ and $v$ has at least two neighbors in $S$ in
$\overline{G^*}$, and the sum of these two numbers is four.
Thus, each has exactly two such neighbors.
By $e(\overline{G^*})=3$, we have
$E(\overline{G^*})=\{uv,ua,vb\}$ for some $a,b\in S\setminus T$.
If $a=b$, then $|\Lambda_{G^*}(a)\cap S|=|S|-1-2=k-1$,
which implies that $a\in T$, contrary to $a\in S\setminus T$.
Hence $a\ne b$, and so $F\cong P_4$.
By $V(F)\subseteq S$ and $|V(F)|=4$, we have $k\ge2$.
According to $t^*=2$, we have $n\le r(k+2)+1$.
Therefore, $F$ satisfies \textup{(f)}.%
\end{subcase}

\begin{subcase}{$(p,t^*)=(3,3)$.}\label{subcase:edge-triangle}

In this case, $t^*x=6=2e(\overline{G^*})$.
By Lemma~\ref{edge:equality} \textup{(ii)}, we have $T\subseteq S$,
and $F=\overline{G^*}[T]$ is a $2$-regular graph of order three.
Hence $F\cong K_3$.
By $t^*=3$, we have $n\le r(k+2)+2$.
Thus, $F$ satisfies \textup{(g)}.%
\end{subcase}

\begin{subcase}{$(p,t^*)=(4,4)$.}\label{subcase:edge-cycle}

In this case, $t^*x=8=2e(\overline{G^*})$.
By Lemma~\ref{edge:equality} \textup{(ii)}, we have $T\subseteq S$,
and $F=\overline{G^*}[T]$ is a $2$-regular graph of order four.
Hence $F\cong C_4$.
Since $V(F)=T\subseteq S$, $4\le k+2$. Hence $k\ge2$.
By $t^*=4$, we have $n\le r(k+2)+3$.
Thus, $F$ satisfies \textup{(h)}.%
\end{subcase}%
\end{case}

Combining Cases \ref{case:edge-zero}-\ref{case:edge-large} with the cases $p=1$ and $x=p$,
we obtain $G^*=K_n-F$ for some $F\in\mathcal F_{n,k,r}$.
Conversely, let $F\in\mathcal F_{n,k,r}$ and put $J=K_n-F$.
Since $e(F)=p\le n-2$, we have $e(J)=\binom n2-p$ and
$\delta(J)\ge n-1-p\ge1$.
It remains to show that $\beta^k(J)<r$.
We verify this for the following four groups of conditions.

\vspace{1mm}
\noindent\textup{(a)} We have $F\cong K_{1,p}$, and hence $J\cong H_{n,p}$.
By Lemma~\ref{lem:p}, we obtain $\beta^k(J)<r$.

\vspace{1mm}
\noindent\textup{(b)} We have $n=(r+1)k+p-1$ and $F\in\mathcal S_{k,p}$.
Let $c\le k$ be the number of components of $F$.
Since $|V(F)|=p+c$,
$|V(J)\setminus V(F)|=n-p-c=(r+1)k-1-c\ge k-c$.
Choose one center from each star and $k-c$ vertices from
$V(J)\setminus V(F)$ to form a set $S$ of size $k$.
The vertices outside $\Lambda_J^k(S)$ are exactly those of $S$
together with the $p$ leaves of these stars.
Hence $\Lambda_J^k(S)\ne V(J)$ and
$\beta^k(J)\le\frac{n-k-p}{k}=\frac{rk-1}{k}<r$.

\vspace{1mm}
\noindent\textup{(c)-(e)} We have $|V(F)|\le k+1\le n$ and
$n\le r(k+1)+|V(F)|-1$.
Choose a vertex set $S$ with $|S|=k+1$ and $V(F)\subseteq S$.
Every vertex of $F$ has at most $k-1$ neighbors in $S$ in $J$,
whereas every vertex outside $V(F)$ has at least $k$ neighbors in $S$.
Thus, $\Lambda_J^k(S)=V(J)\setminus V(F)\ne V(J)$.
By the bound on $n$, we obtain
$\beta^k(J)\le\frac{n-|V(F)|}{k+1}\le r-\frac{1}{k+1}<r$.

\vspace{1mm}
\noindent\textup{(f)-(h)} We have $F\cong P_4$, $K_3$ or $C_4$,
respectively, and $|V(F)|\le k+2$.
If $k\ge2$, then $n\ge2k\ge k+2$.
If $k=1$, then $r\ge2$ and $n\ge r+1\ge3=k+2$.
Hence we may choose a vertex set $S$ with $|S|=k+2$
and $V(F)\subseteq S$.
For every $v\in V(F)$, we have
$|\Lambda_J(v)\cap S|=k+1-d_F(v)$,
while every vertex outside $V(F)$ has at least $k+1$ neighbors in $S$.
Therefore,
$
V(J)\setminus\Lambda_J^k(S)=\{v\in V(F):d_F(v)\ge2\}.
$
The set on the right has $2$, $3$ and $4$ vertices when
$F\cong P_4$, $K_3$ and $C_4$, respectively.
In particular, $\Lambda_J^k(S)\ne V(J)$.
By the bounds on $n$ in \textup{(f)-(h)}, we obtain that
$\beta^k(J)\le\frac{|\Lambda_J^k(S)|}{k+2}
\le r-\frac{1}{k+2}<r$.

Thus, $J\in\mathcal G_{n,k,r}$ and $e(J)=\binom n2-p$
for every $F\in\mathcal F_{n,k,r}$.
Since $G^*$ was chosen arbitrarily among the graphs of maximum size,
we conclude that
$$
e(G)\le e(G^*)=\binom n2-p
$$
with equality if and only if $G\cong K_n-F$ for some
$F\in\mathcal F_{n,k,r}$.
This completes the proof of Theorem~\ref{edge:main}.
\end{proof}

\Needspace{5\baselineskip}

\begin{proof}[\normalfont\bfseries Proof of Theorem~\ref{edge:first}]
Let $G^*$ be a graph with maximum size in $\mathcal G_{n,1,1}$,
and put $q=\floor{\frac{n-1}{2}}$.
By Lemma~\ref{lem:family}, $G^*$ is isomorphic to a spanning subgraph
of $F_{n,t}$ for some integer $t$ with $1\le t\le q$.
By $F_{n,t}\in\mathcal G_{n,1,1}$, we have
$e(G^*)\le e(F_{n,t})\le e(G^*)$.
Thus, $G^*\cong F_{n,t}$.
Put $f(t)=e(F_{n,t})=\binom{n-t-1}{2}+t(t+1)$.
If $n=3$ or $4$, then $q=1$ and
$G^*\cong F_{n,1}=K_1\vee(K_{n-3}\cup2K_1)$.
Hence Theorem~\ref{edge:first} \textup{(iv)} holds for $n\in\{3,4\}$.
In the following, we assume that $n\ge5$. Then $q\ge2$.
Obviously,
$f(t)<\max\{f(1),f(q)\}$ for $1<t<q$.
By the maximality of $G^*$, we obtain $G^*\cong F_{n,1}$ or $F_{n,q}$.
Note that $F_{n,1}=K_1\vee(K_{n-3}\cup2K_1)$.
By $n-2q-1\in\{0,1\}$, we have $F_{n,q}\cong K_q\vee(n-q)K_1$.
Thus, $f(1)=\binom{n-2}{2}+2$ and $f(q)=\binom q2+q(n-q)$.
By a direct calculation, we obtain
$$
f(q)-f(1)=\frac{(q-1)(3q+8-2n)}{2}
=\begin{cases}
\dfrac{(q-1)(6-q)}{2},&n=2q+1,\\[2mm]
\dfrac{(q-1)(4-q)}{2},&n=2q+2.
\end{cases}
$$
We consider the following three cases.

\resetcases
\begin{case}{$n\in\{5,6,7,8,9,11\}$.}\label{case:first-small}

In this case, $f(q)>f(1)$.
Hence $G^*\cong K_q\vee(n-q)K_1$ and
$e(G^*)=\binom q2+q(n-q)$.
This proves Theorem~\ref{edge:first} \textup{(i)}.%
\end{case}

\begin{case}{$n\in\{10,13\}$.}\label{case:first-ties}

If $n=10$, then $q=4$ and $f(1)=f(q)=30$.
The two possible extremal graphs are
$F_{10,1}=K_1\vee(K_7\cup2K_1)$ and $F_{10,4}\cong K_4\vee6K_1$.
If $n=13$, then $q=6$ and $f(1)=f(q)=57$.
The two possible extremal graphs are
$F_{13,1}=K_1\vee(K_{10}\cup2K_1)$ and $F_{13,6}\cong K_6\vee7K_1$.
In either case, both graphs belong to $\mathcal G_{n,1,1}$
and attain the maximum size.
Moreover, they are nonisomorphic, since their minimum degrees
are $1$ and $q$, respectively.
Thus, Theorem~\ref{edge:first} \textup{(ii)} and \textup{(iii)} hold.%
\end{case}

\begin{case}{$n=12$ or $n\ge14$.}\label{case:first-large}

In this case, $f(1)>f(q)$.
Hence $G^*\cong K_1\vee(K_{n-3}\cup2K_1)$ and
$e(G^*)=\binom{n-2}{2}+2$.
Together with the cases $n\in\{3,4\}$, we prove Theorem~\ref{edge:first} \textup{(iv)}.%
\end{case}
\end{proof}

\section{Proofs of Theorems~\ref{thm:main} and \ref{thm:first}}\label{sec:mainproof}

We first present several key lemmas for the proofs of our spectral results.
\begin{lem}[Godsil and Royle \cite{Godsil2001}]\label{lem:pf}
Let $H$ be a spanning subgraph of a connected graph $G$.
Then $\rho(H)\le\rho(G)$ with equality if and only if $H=G$.
\end{lem}

Let $A$ be a square matrix whose rows and columns are partitioned
according to the same partition of the index set.
Write $A_{ij}$ for the resulting blocks.
If every block $A_{ij}$ has constant row sum $q_{ij}$, then the partition
is called equitable, and $Q=(q_{ij})$ is its quotient matrix.

\begin{lem}[Godsil and Royle \cite{Godsil2001}]\label{lem:quotient}
Let $A$ be a nonnegative symmetric matrix with an equitable partition,
and let $Q$ be the corresponding quotient matrix.
Then every eigenvalue of $Q$ is an eigenvalue of $A$.
Moreover, if $A$ is irreducible, then $\rho(A)=\rho(Q)$.
\end{lem}

Recall that $F_{n,t}=K_t\vee\bigl(K_{n-2t-1}\cup(t+1)K_1\bigr)$
for integers $n\ge3$ and $t$ with $1\le t\le\floor{\frac{n-1}{2}}$.
For each integer $n\ge3$, let
$$
H_n=
\begin{cases}
K_2\vee3K_1,&n=5,\\[1mm]
K_2\vee4K_1,&n=6,\\[1mm]
K_3\vee4K_1,&n=7,\\[1mm]
K_4\vee5K_1,&n=9,\\[1mm]
K_1\vee(K_{n-3}\cup2K_1),&n\in\{3,4,8\}\text{ or }n\ge10.
\end{cases}
$$

\begin{lem}\label{lem:Fcompare}
Let $n,t$ be integers with $n\ge3$ and
$1\le t\le\floor{\frac{n-1}{2}}$.
Then $\rho(F_{n,t})\le\rho(H_n)$ with equality if and only if
$F_{n,t}\cong H_n$.
\end{lem}
\begin{proof}
If $n=3$ or $4$, then $t=1$ and $F_{n,t}=H_n$.
Next we consider $n\ge5$.
Put $q=\floor{\frac{n-1}{2}}$.
For $1\le t\le q$, let
$P_{n,t}(z)=z^3-(n-t-3)z^2-(n+t^2-2)z+t(t+1)(n-2t-2)$.
If $n-2t-1>0$, then the vertex sets of $K_t$, $K_{n-2t-1}$
and $(t+1)K_1$, in this order, form an equitable partition of $F_{n,t}$.
The corresponding quotient matrix is
$$
B_{n,t}=\left[\begin{array}{ccc}
t-1&n-2t-1&t+1\\
t&n-2t-2&0\\
t&0&0
\end{array}\right].
$$
By computations, the characteristic polynomial of $B_{n,t}$
is $P_{n,t}(z)$.
Since $F_{n,t}$ is connected, we have
$\rho(F_{n,t})=\rho(B_{n,t})$ by Lemma~\ref{lem:quotient}.
Let $D=\operatorname{diag}(\sqrt t,\sqrt{n-2t-1},\sqrt{t+1})$.
Then $B_{n,t}$ is similar to the real symmetric matrix
$$
DB_{n,t}D^{-1}
=\left[\begin{array}{ccc}
t-1&\sqrt{t(n-2t-1)}&\sqrt{t(t+1)}\\
\sqrt{t(n-2t-1)}&n-2t-2&0\\
\sqrt{t(t+1)}&0&0
\end{array}\right].
$$
Thus, $P_{n,t}$ has three real roots, and its largest root is $\rho(F_{n,t})$.
By deleting the first row and column of $DB_{n,t}D^{-1}$,
we obtain the principal submatrix $\operatorname{diag}(n-2t-2,0)$.
By the eigenvalue interlacing theorem, the second largest root of
$P_{n,t}$ is at most $\max\{n-2t-2,0\}\le n-3$.

If $n-2t-1=0$, then $F_{n,t}=K_t\vee(t+1)K_1$.
The corresponding quotient matrix is
$$
\left[\begin{array}{cc}
t-1&t+1\\
t&0
\end{array}\right].
$$
By Lemma~\ref{lem:quotient}, $\rho(F_{n,t})$ is the positive root of
$z^2-(t-1)z-t(t+1)$.
Note that $P_{n,t}(z)=(z+1)(z^2-(t-1)z-t(t+1))$.
Since $t(t+1)>0$, the two roots of the quadratic factor have opposite signs.
Hence the other two roots of $P_{n,t}$ are real and negative.
Combining the above two cases, we have
$P_{n,t}(z)=(z-\rho(F_{n,t}))(z-\mu_2)(z-\mu_3)$,
where $\mu_3\le\mu_2\le n-3$.
Note that $z-\mu_2>0$ and $z-\mu_3>0$ for $z>n-3$.
Thus, for every real $z>n-3$, we have $P_{n,t}(z)>0$ if and only if
$z>\rho(F_{n,t}).$

\resetcases
\begin{case}{$n=8$ or $n\ge10$.}\label{case:spectral-large}

In this case, $H_n=F_{n,1}$.
If $t=1$, then $F_{n,t}=H_n$.
Next we consider $2\le t\le q$.
Let $\lambda=\rho(F_{n,1})$.
Note that $K_{n-2}\cup2K_1$ is a proper spanning subgraph of
the connected graph $F_{n,1}$.
By Lemma~\ref{lem:pf}, we have $\lambda>n-3$.
Let $g_t(z)=z^2-(t+1)z+(t+2)n-2t^2-6t-8$.
By $P_{n,1}(\lambda)=0$, we obtain
$$
P_{n,t}(\lambda)
=P_{n,t}(\lambda)-P_{n,1}(\lambda)=(t-1)g_t(\lambda).
$$
We show that $g_t(\lambda)>0$.
By $t\le\frac{n-1}{2}$, we have
$g_t'(z)=2z-t-1\ge\frac{3n-13}{2}>0$ for $z\ge n-3$.
This implies that $g_t$ is strictly increasing on $[n-3,\infty)$.
Since $\lambda>n-3$, we obtain $g_t(\lambda)>g_t(n-3)$.
By $n-2t-1\ge0$, we have
\begin{equation}\label{eq:Rbound}
\begin{aligned}
g_t(n-3)
=n^2-5n+4-2t^2-3t
&=\frac{(n-1)(n-10)}{2}
+\frac{(n-2t-1)(n+2t+2)}{2}\\
&\ge\frac{(n-1)(n-10)}{2}.
\end{aligned}
\end{equation}
If $n\ge10$, then $g_t(n-3)\ge0$ by (\ref{eq:Rbound}).
If $n=8$, then $t\in\{2,3\}$, and we have
$g_2(5)=14>0$ and $g_3(5)=1>0$.
Thus, $g_t(\lambda)>0$ in either case.
By $t\ge2$ and $P_{n,t}(\lambda)=(t-1)g_t(\lambda)$, we obtain
$P_{n,t}(\lambda)>0$.
Combining this with $\lambda>n-3$ and the above equivalence, we have
$\rho(F_{n,t})<\lambda=\rho(H_n)$.
Therefore, $\rho(F_{n,t})\le\rho(H_n)$ with equality if and only if $t=1$.%
\end{case}

\begin{case}{$n\in\{5,6,7,9\}$.}\label{case:spectral-small}

By $n-2q-1\in\{0,1\}$, we have
$F_{n,q}\cong K_q\vee(n-q)K_1=H_n$.
If $t=q$, then $\rho(F_{n,t})=\rho(H_n)$.
Next we consider $1\le t<q$.
Let $\lambda_0=\rho(F_{n,q})$.
The vertex sets of $K_q$ and $(n-q)K_1$ form an equitable partition
of $F_{n,q}$.
The corresponding quotient matrix is
$$
\left[\begin{array}{cc}
q-1&n-q\\
q&0
\end{array}\right].
$$
By Lemma~\ref{lem:quotient}, $\lambda_0$ is the positive root of
$z^2-(q-1)z-q(n-q)$.
We show that $\lambda_0>n-3$ and $P_{n,t}(\lambda_0)>0$
for every $1\le t<q$.

If $n=5$, then $q=2$ and $\lambda_0=3>2=n-3$.
By $\lambda_0^2=\lambda_0+6$ and the definition of $P_{n,t}$, we have
$P_{5,1}(\lambda_0)=2\lambda_0+2>0$.
If $n=6$, then $q=2$ and $\lambda_0=\frac{1+\sqrt{33}}{2}>3=n-3$.
By $\lambda_0^2=\lambda_0+8$, we have
$P_{6,1}(\lambda_0)=2\lambda_0-4>0$.
If $n=7$, then $q=3$ and $\lambda_0=1+\sqrt{13}>4=n-3$.
By $\lambda_0^2=2\lambda_0+12$, we obtain
$P_{7,1}(\lambda_0)=4\lambda_0-6>0$ and
$P_{7,2}(\lambda_0)=3\lambda_0+6>0$.
If $n=9$, then $q=4$ and $\lambda_0=\frac{3+\sqrt{89}}{2}>6=n-3$.
By $\lambda_0^2=3\lambda_0+20$, we obtain
$P_{9,1}(\lambda_0)=6\lambda_0-30>0$,
$P_{9,2}(\lambda_0)=6\lambda_0-2>0$ and
$P_{9,3}(\lambda_0)=4\lambda_0+12>0$.

By the above calculations, we have $\lambda_0>n-3$ and
$P_{n,t}(\lambda_0)>0$ for every $1\le t<q$.
Combining this with the above equivalence, we obtain
$\rho(F_{n,t})<\lambda_0=\rho(H_n)$.
Therefore, $\rho(F_{n,t})\le\rho(H_n)$ with equality if and only if $t=q$.%
\end{case}
Combining Cases \ref{case:spectral-large} and \ref{case:spectral-small} with $n\in\{3,4\}$, we have
$\rho(F_{n,t})\le\rho(H_n)$ with equality if and only if
$F_{n,t}\cong H_n$.
This completes the proof.
\end{proof}

For $1\le t\le n-2$, let $\theta_{n,t}=\rho(H_{n,t})$.

\begin{lem}[Rowlinson \cite{Rowlinson1988}]\label{lem:rowlinson}
Let $n,t$ be integers with $1\le t\le n-2$, and let $G$ be a graph
of order $n$ with $e(G)=\binom n2-t$.
Then $\rho(G)\le\theta_{n,t}$ with equality if and only if
$G\cong H_{n,t}$.
\end{lem}

Since $1\le t\le n-2$, the graph $H_{n,t}$ is connected.
Moreover, $K_{n-1}\cup K_1$ is a proper spanning subgraph of $H_{n,t}$,
and $H_{n,t}$ is a proper spanning subgraph of $K_n$.
By Lemma~\ref{lem:pf}, we have $n-2<\theta_{n,t}<n-1.$
We are now ready to prove Theorem~\ref{thm:main}.

\begin{proof}[\normalfont\bfseries Proof of Theorem~\ref{thm:main}]
Let $G^*$ be a graph with maximum spectral radius in $\mathcal G_{n,k,r}$.
Since $G^*\in\mathcal G_{n,k,r}$, we have
$\beta^k(G^*)<r\le\beta^k(K_n)$, and hence $e(\overline{G^*})\ge1$.
By Lemma~\ref{lem:p}, we have $1\le p\le n-2$ and
$\beta^k(H_{n,p})<r$.
Combining this with $\delta(H_{n,p})=n-p-1\ge1$, we obtain
$H_{n,p}\in\mathcal G_{n,k,r}$.
Thus, by the maximality of $G^*$ and $n-2<\theta_{n,t}<n-1,$ we have
\begin{equation}\label{eq:lower}
\rho(G^*)\ge\rho(H_{n,p})=\theta_{n,p}>n-2.
\end{equation}

\begin{claim}\label{clm:spectral-completion}
If $e(\overline{G^*})\ge p$, then $G^*\cong H_{n,p}$.
\end{claim}
\begin{proof}
Let $G'$ be the graph obtained from $G^*$ by adding
$e(\overline{G^*})-p$ edges.
Then $e(G')=\binom n2-p>\binom{n-1}{2}$.
Since every disconnected graph of order $n$ has at most
$\binom{n-1}{2}$ edges, $G'$ is connected.
By Lemmas~\ref{lem:pf} and \ref{lem:rowlinson}, we have
$\rho(G^*)\le\rho(G')\le\theta_{n,p}$.
Combining this with (\ref{eq:lower}), we obtain
$\rho(G^*)=\rho(G')=\theta_{n,p}$.
Note that $G^*$ is a spanning subgraph of the connected graph $G'$
and $\rho(G^*)=\rho(G').$ Then, by
Lemma~\ref{lem:pf}, we have $G^*=G'$.
Moreover, $e(G')=\binom n2-p$ and
$\rho(G')=\theta_{n,p}=\rho(H_{n,p})$.
By Lemma~\ref{lem:rowlinson},
we obtain that $G'\cong H_{n,p}$.
Thus, $G^*\cong H_{n,p}$.
\end{proof}

\vspace{1mm}
\noindent\textup{(i)} Note that $(n,k,r)=(11,3,2)$.
By the definition of $p$, we have $p=3$.
Let $G'=K_{11}-2K_2$, and let $S_0$ be the set of four vertices
incident with the deleted edges.
As shown in the proof of Theorem~\ref{edge:main}, we have
$G'\in\mathcal G_{11,3,2}$.
We show that $\rho(G')>\theta_{11,3}$.
Every vertex of $S_0$ has two neighbors in $S_0$ and seven outside $S_0$,
while every vertex outside $S_0$ has four neighbors in $S_0$
and six outside $S_0$.
Hence the partition $\{S_0,V(G')\setminus S_0\}$ is equitable.
For $H_{11,3}=K_7\vee(K_3\cup K_1)$, the vertex sets of
$K_7$, $K_3$ and $K_1$, in this order, also form an equitable partition.
The corresponding quotient matrices are
$$
R=\left[\begin{array}{cc}
2&7\\
4&6
\end{array}\right],
\qquad
Q=\left[\begin{array}{ccc}
6&3&1\\
7&2&0\\
7&0&0
\end{array}\right].
$$
By a direct calculation, we have $\det(zI-R)=z^2-8z-16$ and
$\det(zI-Q)=z^3-8z^2-16z+14$.
By Lemma~\ref{lem:quotient}, we obtain
$\lambda=\rho(G')=4+4\sqrt2$, and $\theta_{11,3}$ is the largest root
of $f(z)=z^3-8z^2-16z+14$.
Since $\lambda^2-8\lambda-16=0$, we have $f(\lambda)=14>0$.
Moreover, $f'(z)=3z^2-16z-16\ge83>0$ for $z\ge9$.
By $n-2<\theta_{n,t}<n-1$, we have $\theta_{11,3}>9$.
Combining this with $\lambda>9$ and
$f(\theta_{11,3})=0<f(\lambda)$, we obtain
$$
\rho(G')=\lambda>\theta_{11,3}=\rho(H_{11,3}).
$$
By the maximality of $G^*$, we have
$\rho(G^*)\ge\rho(G')>\theta_{11,3}$.
By Claim~\ref{clm:spectral-completion}, we have
$e(\overline{G^*})<p$.
By Lemmas~\ref{lem:p} and \ref{lem:exception}, we obtain
$G^*\cong K_{11}-2K_2$.
Hence $\rho(G)\le\rho(K_{11}-2K_2)$ with equality if and only if
$G\cong K_{11}-2K_2$.
This proves \textup{(i)}.

\vspace{1mm}
\noindent\textup{(ii)} Note that $(n,k,r)\ne(11,3,2)$.
Suppose to the contrary that $e(\overline{G^*})<p$.
By Lemma~\ref{lem:p}, we have
$\beta^k(G^*)<r\le\beta^k(H_{n,e(\overline{G^*})})$.
By Lemma~\ref{lem:exception}, we obtain $(n,k,r)=(11,3,2)$,
a contradiction.
Hence $e(\overline{G^*})\ge p$.
By Claim~\ref{clm:spectral-completion}, we have $G^*\cong H_{n,p}$.
Since $G^*$ was chosen arbitrarily among the graphs with maximum
spectral radius, we conclude that
$\rho(G)\le\rho(H_{n,p})$ with equality if and only if $G\cong H_{n,p}$.
This proves \textup{(ii)} and completes the proof of Theorem~\ref{thm:main}.
\end{proof}

We next consider Theorem~\ref{thm:first}.

\begin{proof}[\normalfont\bfseries Proof of Theorem~\ref{thm:first}]
Let $G^*$ be a graph with maximum spectral radius in $\mathcal G_{n,1,1}$.
By Lemma~\ref{lem:family}, $G^*$ is isomorphic to a spanning subgraph
of $F_{n,t}$ for some integer $t$ with
$1\le t\le\floor{\frac{n-1}{2}}$.
Since $H_n$ is isomorphic to one of the graphs in this family,
we have $H_n\in\mathcal G_{n,1,1}$ by Lemma~\ref{lem:family} \textup{(ii)}.
By the maximality of $G^*$, we obtain $\rho(G^*)\ge\rho(H_n)$.
Moreover, $F_{n,t}$ is connected.
By Lemmas~\ref{lem:pf} and \ref{lem:Fcompare}, we have
$$
\rho(H_n)\le\rho(G^*)\le\rho(F_{n,t})\le\rho(H_n).
$$
By $\rho(G^*)=\rho(F_{n,t})$ and the equality condition in
Lemma~\ref{lem:pf}, we have $G^*\cong F_{n,t}$.
By $\rho(F_{n,t})=\rho(H_n)$ and the equality condition in
Lemma~\ref{lem:Fcompare}, we obtain $F_{n,t}\cong H_n$.
Thus, $G^*\cong H_n$.
Since $G^*$ was chosen arbitrarily among the graphs with maximum
spectral radius, we conclude that
$\rho(G)\le\rho(H_n)$ with equality if and only if $G\cong H_n$.
This completes the proof of Theorem~\ref{thm:first}.
\end{proof}

\vspace{5mm}
\noindent
{\bf Declaration of competing interest}
\vspace{3mm}

The authors declare that they have no known competing financial interests or personal relationships that could have appeared to influence the work reported in this paper.

\vspace{5mm}
\noindent
{\bf Data availability}
\vspace{3mm}

No data was used for the research described in this paper.

%

\end{document}